\documentclass[sn-mathphys-ay]{sn-jnl}

\usepackage{graphicx}
\usepackage{multirow}
\usepackage{amsmath,amssymb,amsfonts}
\usepackage{amsthm}
\usepackage{mathrsfs}
\usepackage[title]{appendix}
\usepackage{xcolor}
\usepackage{textcomp}
\usepackage{manyfoot}
\usepackage{booktabs}
\usepackage{algorithm}
\usepackage{algorithmicx}
\usepackage{algpseudocode}
\usepackage{listings}
\usepackage{cleveref}

\newtheorem{theorem}{Theorem}

\newtheorem{definition}{Definition}

\begin{document}

\title[Parameter Estimation and Identifiability in Kinetic Flux Profiling Models of Metabolism]{Parameter Estimation and Identifiability in Kinetic Flux Profiling Models of Metabolism}

\author[1,2]{\fnm{Breanna} \sur{Guppy}}\email{breanna-guppy@uiowa.edu}

\author*[1,2]{\fnm{Colleen} \sur{Mitchell}}\email{colleen-mitchell@uiowa.edu}

\author[3,4,5,6]{\fnm{Eric B.} \sur{Taylor}}\email{eric-taylor@uiowa.edu}

\affil*[1]{\orgdiv{Mathematics}, \orgname{University of Iowa}, \orgaddress{\street{2 West Washington Street}, \city{Iowa City}, \postcode{52242}, \state{Iowa}, \country{USA}}}

\affil[2]{\orgdiv{Applied Mathematical \& Computational Sciences}, \orgname{University of Iowa}, \orgaddress{\street{2 West Washington Street}, \city{Iowa City}, \postcode{52242}, \state{Iowa}, \country{USA}}}

\affil[3]{\orgdiv{Molecular Physiology and Biophysics}, \orgname{University of Iowa}, \orgaddress{\street{51 Newton Road}, \city{Iowa City}, \postcode{52242}, \state{Iowa}, \country{USA}}}

\affil[4]{\orgdiv{Fraternal Order of Eagles Diabetes Research Center (FOEDRC)}, \orgname{University of Iowa}, \orgaddress{\street{169 Newton Road}, \city{Iowa City}, \postcode{52246}, \state{Iowa}, \country{USA}}}

\affil[5]{\orgdiv{Pappajohn Biomedical Institute}, \orgname{University of Iowa}, \orgaddress{\street{169 Newton Road}, \city{Iowa City}, \postcode{52246}, \state{Iowa}, \country{USA}}}

\affil[6]{\orgdiv{University of Iowa Carver College of Medicine}, \orgname{University of Iowa}, \orgaddress{\street{451 Newton Road}, \city{Iowa City}, \postcode{52246}, \state{Iowa}, \country{USA}}}

\abstract{Metabolic fluxes are the rates of life-sustaining chemical reactions within a cell and metabolites are the components. 
Determining the changes in these fluxes is crucial to understanding diseases with metabolic causes and consequences. 
Kinetic flux profiling (KFP) is a method for estimating flux that utilizes data from isotope tracing experiments. 
In these experiments, the isotope-labeled nutrient is metabolized through a pathway and integrated into the downstream metabolite pools. 
Measurements of proportion labeled for each metabolite in the pathway are taken at multiple time points and used to fit an ordinary differential equations model with fluxes as parameters.
We begin by generalizing the process of converting diagrams of metabolic pathways into mathematical models composed of differential equations and algebraic constraints. 
The scaled differential equations for proportions of unlabeled metabolite contain parameters related to the metabolic fluxes in the pathway. 
We investigate flux parameter identifiability given data collected only at the steady state of the differential equation. 
Next, we give criteria for valid parameter estimations in the case of a large separation of timescales with fast-slow analysis. 
Bayesian parameter estimation on simulated data from KFP experiments containing both irreversible and reversible reactions illustrates the accuracy and reliability of flux estimations. 
These analyses provide constraints that serve as guidelines for the design of KFP experiments to estimate metabolic fluxes.}


\maketitle

\section{Introduction}\label{sec1}

Metabolism refers to the interconnected set of chemical reactions that allow an organism to utilize nutrients, break down waste, grow, and function. 
The substrates and products of these reactions are called metabolites and the reaction rates are called metabolic fluxes. 
More abstractly, metabolic flux can be viewed as the rate at which metabolites flow between states in a network representation of a metabolic pathway. 
In response to changes in diet, cellular stress, or disease, distinct patterns emerge in metabolite concentrations indicating changes in the metabolic fluxes. 
Therefore, metabolite concentrations and fluxes are essential information for defining cell physiology and investigating changes in metabolic function \citep{emwas2022fluxomics,stephanopoulos1999metabolic}. 
Since disruption in metabolic networks is connected to numerous disease states including cancer, heart disease, and diabetes, the ability to measure changes in metabolic fluxes can provide insight into the pathophysiology of these conditions \citep{hollywood2006metabolomics}. 
Quantification of metabolic fluxes has the potential to improve diagnosis of diseases with metabolic features and aid in the development of novel treatments \citep{Agus2020GutMM, AronWisnewsky2020MetabolismAM,gowda2014overview, Zheng2021HyocholicAS}.   

Several methods for estimating metabolic fluxes have been proposed. 
Earlier methods, referred to as Flux Balance Analysis, were based on constrained optimization and required the assumption that the system was optimizing some element of cellular function or production \citep{orth2010flux}. 
Metabolic Flux Analysis is also based on constrained optimization and incorporates some measurable rates such as inputs, outputs, or oxygen consumption to further constrain the unknown metabolic fluxes \citep{antoniewicz2015methods}. 
Another similar approach, $^{13}$C Metabolic Flux Analysis, \citep{antoniewicz2015methods} incorporates measurement of levels of $^{13}$C, a stable isotope of Carbon. 
Metabolites in which the carbons have been replaced with this heavy isotope are called labeled and the small change in the mass can be detected using mass spectrometry (MS) or nuclear magnetic resonance (NMR). Following the incorporation of a labeled nutrient, proportions of each metabolite labeled are measured after the labeling patterns have equilibrated. This steady-state pattern of isotope incorporation is included as an additional constraint. 

An advance from these methods is Kinetic Flux Profiling (KFP). This method takes advantage of the dynamics of the changing isotope labeling patterns to better determine the network fluxes. 
This combined experimental and computational approach is sometimes referred to as dynamic flux analysis \citep{xiong2020dynamic}. 
The key principle behind KFP is that greater metabolic fluxes are linked to quicker transmission of isotopic label from a labeled nutrient input \citep{yuan2008kfp}. 
An ordinary differential equation (ODE) model for isotopic labeling is derived from the metabolic network architecture. 
By fitting the parameters in that model, we obtain estimates of the underlying fluxes. 
The method requires that the metabolites in the network maintain constant concentration levels, but the isotopically labeled proportion of each metabolite changes over time. 
We therefore say that the network is in a metabolic steady state but not in an isotopic steady state. Some examples of the successful application of KFP include the investigation of changes in metabolic fluxes in E. coli during starvation \citep{yuan2006kinetic, yuan2007differentiating} and in oleaginous green algae \citep{wu2016kinetic}.

Kinetic Flux Profiling is an ideal forum for the exploration of parameter identifiability. 
The method is intended, from its inception, to provide data that can be used to estimate certain parameters in an ODE system. 
These parameters are the central goal of the profiling technique because of the insight they can provide into the function of a metabolic network. 
While the setup of the system of ODEs from a directed weighted graph is not unique to this system, KFP models have several idiomatic characteristics that give the equations additional structure. 
These include the fluxes of both labeled and unlabeled substrates into the pathway, the mixed fluxes out of the pathway, and the addition of algebraic constraints on the flux parameters. 
In the next section we will give a brief introduction to the KFP experimental procedure.

\subsection{Kinetic Flux Profiling}
The method of Kinetic Flux Profiling includes experimental data collection, conversion of a network diagram to a model consisting of a system of ordinary differential equations and algebraic constraints, and a computational step for the estimation of parameters from the data. 
The experiment begins by allowing the cells to reach a metabolic steady state in unlabeled media. 
At the metabolic steady state, all fluxes and all metabolite concentrations are constant. 
Next, the unlabeled media is switched out for its stable-isotope-labeled equivalent. 
In most cases, one specific metabolite in the media will be labeled with the carbon isotope $^{13}$C in place of all its carbons while the rest of the media remains unchanged. 
After the switch to isotope-labeled media, the $^{13}$C carbons will be incorporated into metabolites within the network. 
Samples from the cells are taken at multiple time points, and the proportions of labeled and unlabeled metabolites are quantified using MS or NMR. 
Ideally, these steps are performed without disrupting the metabolic steady state. 
That is, without perturbing the metabolic concentrations or fluxes.

Next, the architecture of the metabolic pathway under study is used to write a system of ordinary differential equations that describe the change in isotope labeling. 
The metabolic pathway is represented as a directed weighted graph where the nodes represent metabolites and the edges represent reactions. 
The pathway also includes additional inputs and outputs, edges that do not begin or do not terminate at a node. 
To ensure non-trivial solutions, the pathway will contain at least one labeled input and at least one output. 
For simplicity, the diagram should only include metabolites that may become labeled. 
That is, there must be a path from a labeled input to each metabolite in the diagram. 
In a metabolic steady state the concentrations of the metabolites are not changing. 
Therefore the total flux into and out of each metabolite must be equal. 
This flux balance condition gives an additional set of algebraic constraints on the fluxes. 

The resulting system contains parameters for each flux in the network as well as for the total concentration of each metabolite. 
In a standard KFP protocol, it is presumed that these total concentrations are all known. 
The final step is to use a parameter fitting method to fit the parameters for the metabolic fluxes to the isotope labeling time courses. 
Since the ODE system is linear, if we had access to perfect, noise-free, continuous time data the parameters could, in theory, be computed. 
However, we address the accuracy and reliability of the parameter estimates when they are made with limited and noisy data. 
In the following sections, we seek to answer three crucial questions about the estimation of fluxes in KFP, motivated by realistic limitations on data availability. 
What can we learn without estimates of the total concentrations? 
What can we learn with only steady-state labeled proportions? 
What can we learn when there is a separation of time scales between metabolites? 

We begin by introducing a general method for formulating the ODE system and flux parameter constraints. 
We next suggest a scaling of the system that does not include explicit dependence on the total concentrations of each metabolite. 
We finish \cref{sec2} with a condition on the graph that allows us to estimate relative fluxes (proportion of each flux to the total through a metabolite) and turnover rates of each metabolite without measuring metabolite concentrations. 
In \cref{sec3}, we investigate the steady-state problem. 
Often metabolic changes occur on a time scale that make it technically challenging to collect sufficient data before isotopic steady state is reached. 
While it is not possible, in general, to compute turnover rates from steady-state data, we give a simple necessary condition for recovering the relative fluxes from steady-state only. 
This is illustrated with two examples using Bayesian parameter estimation on simulated data to show both the accuracy and reliability of the estimates. 
Finally, in \cref{sec4} we address what can happen if one of the turnover rates is substantially faster than the others. 
We show that the rapid turnover rate parameter is not identifiable unless the fast metabolite is the one that receives labeled input. 
Again, this is illustrated using Bayesian parameter estimation on simulated data from specific examples.  

\section{The KFP Model}\label{sec2}

\begin{figure}[H]
	\centering
	\includegraphics[width=0.75\textwidth]{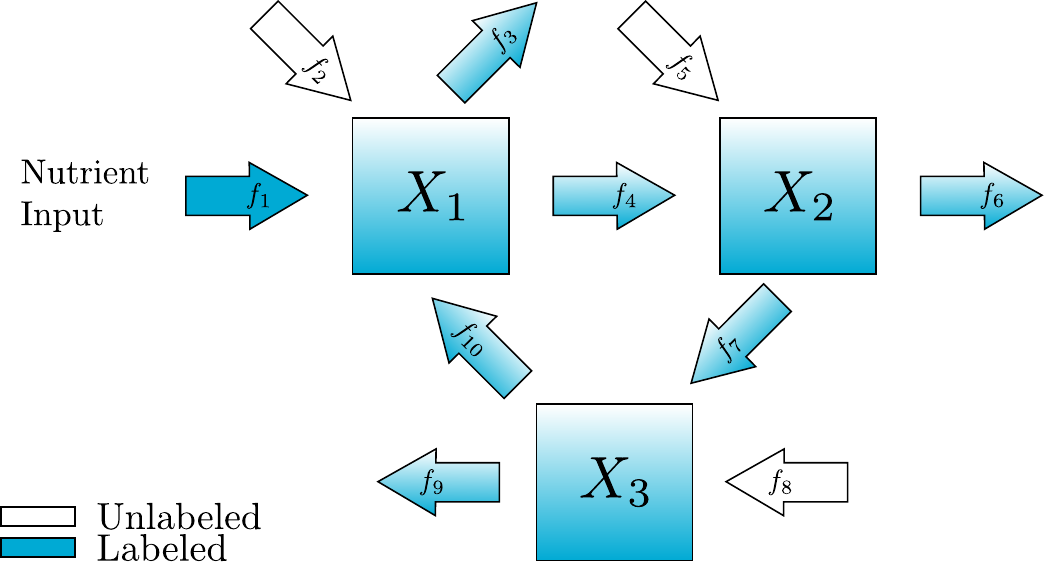}
	\caption{An example reaction network diagram with three metabolites connected in a cycle. There is one influx which is labeled. Each metabolite also has an unlabeled influx and a partially labeled outflux}
    \label{Cyclic}
\end{figure}

In this section, we introduce a new standard procedure for converting the diagram representation of a metabolic network into a matrix representation of the differential equation and algebraic constraints. 
For clarity, we begin with a simple example depicted in \cref{Cyclic} and written in the original notation typically used to describe KFP equations \citep{yuan2008kfp}. 
In this example there are three metabolites, named $X_1$, $X_2$, and $X_3$ and ten reactions with constant fluxes $f_1, f_2, \dots f_{10}$. 
We can write a system of differential equation in terms of the concentrations of either labeled metabolites $x_i^{\mathcal{L}}$ or unlabeled metabolites $x_i^{\mathcal{U}}$, $i=1,2,3$. 
Here we follow \citep{yuan2008kfp} and write all equations in terms of the unlabeled metabolites and use the notation $x_i^{\mathcal{T}}=x_i^{\mathcal{U}}+x_i^{\mathcal{L}}$ for the total concentration of each metabolite. 
In this context, $x_i^{\mathcal{L}}$, $x_i^{\mathcal{U}}$ and $x_i^{\mathcal{T}}$ are all in units of concentration and the fluxes are in units of concentration over time. 

Next, we write a system of differential equations consisting of terms for each flux into or out of each metabolite. 
Terms that correspond to fluxes originating from metabolite $i$ will contain unlabeled metabolite proportional to the current proportion unlabeled for that metabolite $x_i^{\mathcal{U}}/x_i^{\mathcal{T}}$. 
Terms corresponding to fluxes that do not originate at a metabolite are entirely unlabeled. 

    \begin{equation}\label{ExampleDE}
        \begin{aligned}
			x^{\mathcal{U}}_1 &=  
            \frac{x_3^{\mathcal{U}}}{x_3^{\mathcal{T}}} f_{10} - 
            \frac{x_1^{\mathcal{U}}}{x_1^{\mathcal{T}}} f_{4} -
            \frac{x_1^{\mathcal{U}}}{x_1^{\mathcal{T}}} f_{3} + f_2  \\  
            x^{\mathcal{U}}_2 &=  
            \frac{x_1^{\mathcal{U}}}{x_1^{\mathcal{T}}} f_{4} - 
            \frac{x_2^{\mathcal{U}}}{x_2^{\mathcal{T}}} f_{7} -
            \frac{x_2^{\mathcal{U}}}{x_2^{\mathcal{T}}} f_{6} + f_5\\
            x^{\mathcal{U}}_3 &=  
            \frac{x_2^{\mathcal{U}}}{x_2^{\mathcal{T}}} f_{7} - 
            \frac{x_3^{\mathcal{U}}}{x_3^{\mathcal{T}}} f_{10} -
            \frac{x_3^{\mathcal{U}}}{x_3^{\mathcal{T}}} f_{9} + f_8  
		\end{aligned}
    \end{equation}
    
The initial conditions are that all metabolites begin as unlabeled 
    \begin{equation}\label{ExampleIC}
		x_i^{\mathcal{U}}(0)=x_i^{\mathcal{T}} \;\; i=1,2,3.
    \end{equation}
Note that these equations do not contain any term for the labeled input $f_1$ because this influx contributes only to the labeled pool. 
Therefore, this flux would appear in the system for labeled concentrations. It is essential because without it, no label will be added and all metabolites will remain unlabeled. 
A similar formulation is possible for the labeled concentrations and combined for the total concentrations. 
For existing KFP methods, the additional assumption is made that the total concentrations of each metabolite are not changing. 
This assumption means that the fluxes into and out of each metabolite must balance leading to a set of constraints on the fluxes. 
In this example 
    \begin{equation}\label{ExampleConstraints}
        \begin{aligned}
	        f_1+f_2+f_{10} &=f_4 +f_3\\
            f_4+f_5 &=f_6+f_7\\
            f_7+f_8 &=f_9+f_{10}
		\end{aligned}
    \end{equation}

While classic KFP analysis requires that the total concentrations be constant and thus induces these algebraic constraints, the formulation of the equations \cref{ExampleDE} can be generalized to a case where the fluxes are constant but total concentrations $x_i^{\mathcal{T}}$ are allowed to change. 
In this case, each $x_i^{\mathcal{T}}$ will be a linear function of time and the equations \cref{ExampleDE}, while still linear, will not be autonomous. 

The purpose of the remainder of this section is to formalize this process of transforming the graph representation of the reaction diagram into a system of equations in the form of equations \cref{ExampleDE}, \cref{ExampleIC}, \cref{ExampleConstraints}. The formulation of the ODE system in existing literature is ad hoc for each individual application. In these applications, it is typically assumed that the network will be an arborescence with one labeled input and exactly one directed walk to each node. Here we make a more general formulation that can also be used when the graph contains cycles such as in our example or reversible reactions as in the examples in \cref{sec3} and \cref{sec4}.
Our goal is to write the ODE system for isotope-labeling change for a general diagram of a pathway containing \(N\) metabolites and $R$ reactions in the form
	\begin{equation}\label{ODE}
		\dot{\boldsymbol{X}}^\mathcal{U}= \boldsymbol{\hat{A}}\boldsymbol{X}^\mathcal{U}+\boldsymbol{\hat{b}}.
	\end{equation}
In \cref{ODE} the vector \(\boldsymbol{X}^\mathcal{U} \in \mathbb{R}^{N}\) represents the concentrations of unlabeled metabolites at time \(t\). 
The matrix \(\boldsymbol{\hat{A}}\in \mathbb{R}_{\geq 0}^{N\times N}\) describes the biochemical reactions in the pathway whose fluxes originate from a metabolite. This includes the three types of fluxes that will be addressed individually before combining into matrix $\boldsymbol{\hat{A}}$. 
There are positive terms that are inputs from other metabolites, negative terms that are outputs to other metabolites, and negative terms that correspond to fluxes out of the system. 
The vector \(\boldsymbol{\hat{b}} \in \mathbb{R}^N\) consists of the fluxes providing unlabeled metabolite into the pathway from a source outside the scope of the diagram. 

We denote the time of the media switch from unlabeled to labeled as time zero. 
Therefore the initial conditions are that all metabolites begin as unlabeled. 
\begin{equation*}
			x_i^{\mathcal{U}}(0)=x_i^{\mathcal{T}} \;\; i=1, \cdots N
   \end{equation*}

Additionally, we derive the algebraic constraints for metabolic steady state. 
The algebraic constraints are defined by the matrix equation
	\begin{equation}\label{constraints}
		\boldsymbol{F}_{\text{in}} = \boldsymbol{F}_{\text{out}} = \boldsymbol{F}.
	\end{equation}
The diagonal matrices \(\boldsymbol{F}_{\text{in}}, \boldsymbol{F}_{\text{out}} \in \mathbb{R}^{N\times N} \) consist of the total flux into and out of each metabolite. 
When we additionally assume that these fluxes are equal, we can introduce the diagonal matrix $\boldsymbol{F}$ for the flux through each metabolite. 
Descriptions of $\boldsymbol{\hat{A}}$, $\boldsymbol{\hat{b}}$,  \(\boldsymbol{F}_\text{in}\), \(\boldsymbol{F}_\text{out}\), and $\boldsymbol{F}$ based on the structure of the metabolic pathway are given in the following section. 

\subsection{Graph representation of the pathway diagram}

Suppose the diagram of the metabolic pathway is represented by the graph \(G\). 
Let \(G = (\mathcal{N}, \mathcal{E})\) be a special directed weighted graph where the set of nodes \(\mathcal{N}\) represent metabolites and the set of edges \(\mathcal{E}\) have weights that represent metabolic flux in the direction of the edge. 
The set \(\mathcal{E}\) contains edges that do not originate or terminate at a node (i.e. they enter and exit the scope of the graph). It is therefore convenient to partition \(\mathcal{E}\) into four disjoint subsets
    \[\mathcal{E} = \mathcal{E}_\mathcal{L} \cup \mathcal{E}_\mathcal{U} \cup \mathcal{E}_\mathcal{V} \cup \mathcal{E}_\mathcal{W} \]
where \(\mathcal{E}_\mathcal{L}\) is the set of edges that provide labeled metabolite to the pathway, \(\mathcal{E}_\mathcal{U}\) is the set of edges that provide unlabeled metabolite to the pathway, \(\mathcal{E}_\mathcal{V}\) is the set of edges that exit the scope of the pathway, and \(\mathcal{E}_\mathcal{W}\) are the edges that connect nodes within the pathway. 
We will consider only connected graphs in which there is a path from the labeled input(s) to every node. This is because any node that can not receive label will not provide any additional information. 
In some cases where a system has edges exiting or entering, another node is added to be the target or origin of those edges. 
We have chosen not to do that in this case both because it would require separate handling of inputs with labeled and unlabeled metabolites and because we feel the notion of entering and exiting the pathway is useful for the interpretation of the experiments.

We will use the following notation for the number of metabolites and reactions in the pathway, \(|\mathcal{N}| = N\) and \(|\mathcal{E}| = R\). 
Because we only consider the case where there is at least one labeled influx $1\le|\mathcal{E}_{\mathcal{L}}|\le N$. 
Each metabolite may have an unlabeled influx and there must be at least one outflux so $0\le |\mathcal{E}_{\mathcal{U}}| \le N$ and $1\le |\mathcal{E}_{\mathcal{V}}|\le N$. 
Finally, since the graph is connected, $N-1 \le |\mathcal{E}_{\mathcal{W}}|\le N(N-1)$. 
In Figure \ref{Cyclic} we have $N=3$ metabolites and $R=10$ reactions. 
$\mathcal{E}_{\mathcal{L}}$ corresponds to the labeled inflow with rate $f_1$. 
$\mathcal{E}_{\mathcal{U}}$ corresponds to the unlabeled inflows with rates $f_2, f_5$ and $f_8$. 
$\mathcal{E}_{\mathcal{V}}$ corresponds to the outflows with rates $f_3, f_6$ and $f_9$. 
$\mathcal{E}_{\mathcal{W}}$ corresponds to the flows between nodes with rates $f_4, f_7$ and $f_{10}$.

By considering only the directed weighted graph with edges in $\mathcal{E}_{\mathcal{W}}$ we can define the weighted adjacency matrix.
        The \textbf{weighted adjacency matrix}, \(\boldsymbol{W}\in\mathbb{R}^{N\times N}_{\geq 0}\), is the matrix composed of the weights of the directed edges between nodes, $\mathcal{E}_{\mathcal{W}}$. The matrix entry \(w_{i,j}\) is the weight of the edge from node \(i\) to node \(j\). If no connection between node \(i\) and node \(j\) exists \(w_{ij}=0\) \citep{musulin2014process}. 
We further define the total concentration matrix $\boldsymbol{X}^{\mathcal{T}}$ as the diagonal matrix containing the total concentration of each metabolite ${[\boldsymbol{X}^{\mathcal{T}}]}_{i,i}=x_i^{\mathcal{T}}$. 
Then the terms in \cref{ODE} corresponding to the influx into each node from the other nodes is $\boldsymbol{W}^T {\boldsymbol{X}^{\mathcal{T}}}^{-1} \boldsymbol{X}^{\mathcal{U}}$. In our example,

    \begin{equation*}
        \boldsymbol{W}^T {\boldsymbol{X}^{\mathcal{T}}}^{-1} \boldsymbol{X}^{\mathcal{U}} = \begin{pmatrix}
        0 & 0 & f_{10} \\ f_4 & 0 & 0 \\
            0 & f_7 & 0
        \end{pmatrix}
        \begin{pmatrix}
            \frac{1}{x_1^{\mathcal{T}}} & 0 & 0\\
            0 & \frac{1}{x_2^{\mathcal{T}}} & 0\\
            0 & 0 & \frac{1}{x_3^{\mathcal{T}}}
        \end{pmatrix}
        \begin{pmatrix}
            x_1^{\mathcal{U}} \\x_2^{\mathcal{U}}\\x_3^{\mathcal{U}}
        \end{pmatrix}
        = \begin{pmatrix}
             f_{10} \frac{x_3^{\mathcal{U}}}{x_3^{\mathcal{T}}}\\
           f_4 \frac{x_1^{\mathcal{U}}}{x_1^{\mathcal{T}}}\\
           f_7 \frac{x_2^{\mathcal{U}}}{x_2^{\mathcal{T}}}
        \end{pmatrix}.
    \end{equation*}
    
Another commonly used matrix in graph theory is the degree matrix. 
The {\textit{degree matrix}}, \(\boldsymbol{D}\), is a diagonal matrix of the degree of each node \citep{liproperties}. 
Since we have a directed weighted graph, we will define two variations of the degree matrix: the weighted out-degree matrix and the weighted in-degree matrix.

	\begin{definition} 
		The {\textit{weighted out-degree matrix}}, \(\boldsymbol{D}_\text{out}\in\mathbb{R}^{N\times N}\), is a diagonal matrix composed of the sum of weights on the edges leaving a node and entering another node. The matrix entries are computed by \({[\boldsymbol{D}_\text{out}]}_{i,i}= \sum_{j=1}^n w_{i,j}\).
		Similarly, the {\textit{weighted in-degree matrix}}, \(\boldsymbol{D}_\text{in}\in\mathbb{R}^{N\times N}\), is a diagonal matrix composed of the sum of weights on the edges entering a node from another node. The matrix entries are computed by \({[\boldsymbol{D}_\text{in}]}_{j,j}=\sum_{i=1}^n w_{i,j}\).
	\end{definition}

The terms in \cref{ODE} that correspond to outfluxes to other metabolites can be written. $\boldsymbol{D}_{\text{out}}{\boldsymbol{X}^{\mathcal{T}}}^{-1} \boldsymbol{X}^{\mathcal{U}}$. In our example,

    \begin{equation} 
        \boldsymbol{D}_{\text{out}}{\boldsymbol{X}^{\mathcal{T}}}^{-1} \boldsymbol{X}^{\mathcal{U}}= \begin{pmatrix}
        f_4 & 0 & 0 \\ 0 & f_7 & 0 \\
            0 & 0 & f_{10}
        \end{pmatrix}
        \begin{pmatrix}
            \frac{1}{x_1^{\mathcal{T}}} & 0 & 0\\
            0 & \frac{1}{x_2^{\mathcal{T}}} & 0\\
            0 & 0 & \frac{1}{x_3^{\mathcal{T}}}
        \end{pmatrix}
        \begin{pmatrix}
            x_1^{\mathcal{U}} \\x_2^{\mathcal{U}}\\x_3^{\mathcal{U}}
        \end{pmatrix}
        = \begin{pmatrix}
             f_{4} \frac{x_1^{\mathcal{U}}}{x_1^{\mathcal{T}}}\\
           f_7 \frac{x_2^{\mathcal{U}}}{x_2^{\mathcal{T}}}\\
           f_{10} \frac{x_3^{\mathcal{U}}}{x_3^{\mathcal{T}}}
        \end{pmatrix}.
    \end{equation}

Our graph \(G\) also includes edges in sets \(\mathcal{E}_\mathcal{L}\), \(\mathcal{E}_\mathcal{U}\), and \(\mathcal{E}_\mathcal{V}\) that enter and exit the scope of the graph so we define matrices corresponding to these sets of fluxes.

   \begin{definition}
	   The \textbf{labeled entry matrix}, $\boldsymbol{D}_{\mathcal{L}} \in\mathbb{R}^{N\times N}$, is the diagonal matrix consisting of the weights of the edges in \(\mathcal{E}_\mathcal{L}\), ${[\boldsymbol{D}_\mathcal{L}]}_{i,i}$ is the weight of the edge that enters nodes \(i\). The \textbf{unlabeled entry matrix}, $\boldsymbol{D}_{\mathcal{U}} \in\mathbb{R}^{N\times N}$, is the diagonal matrix consisting of the weights of the edges in \(\mathcal{E}_\mathcal{U}\) for unlabeled metabolites that enter nodes \(i\). 
	   The \textbf{exit matrix},  \(\boldsymbol{D}_{\mathcal{V}}\in\mathbb{R}^{N\times N}\), is the diagonal matrix consisting of the weights in \(\mathcal{E}_\mathcal{V}\) for the edges that leave node \(i\).
    \end{definition}

The terms in \cref{ODE} which correspond to outfluxes exiting the graph can be written $\boldsymbol{D}_\mathcal{V}{\boldsymbol{X}^{\mathcal{T}}}^{-1} \boldsymbol{X}^{\mathcal{U}}$. In our example,

    \begin{equation*} 
        \boldsymbol{D}_{\mathcal{V}}{\boldsymbol{X}^{\mathcal{T}}}^{-1} \boldsymbol{X}^{\mathcal{U}}= \begin{pmatrix}
        f_3 & 0 & 0 \\ 0 & f_6 & 0 \\
            0 & 0 & f_{9}
        \end{pmatrix}
        \begin{pmatrix}
            \frac{1}{x_1^{\mathcal{T}}} & 0 & 0\\
            0 & \frac{1}{x_2^{\mathcal{T}}} & 0\\
            0 & 0 & \frac{1}{x_3^{\mathcal{T}}}
        \end{pmatrix}
        \begin{pmatrix}
            x_1^{\mathcal{U}} \\x_2^{\mathcal{U}}\\x_3^{\mathcal{U}}
        \end{pmatrix}
        = \begin{pmatrix}
             f_{3} \frac{x_1^{\mathcal{U}}}{x_1^{\mathcal{T}}}\\
           f_6 \frac{x_2^{\mathcal{U}}}{x_2^{\mathcal{T}}}\\
           f_{9} \frac{x_3^{\mathcal{U}}}{x_3^{\mathcal{T}}}
        \end{pmatrix}.
    \end{equation*}

The last set of terms in \cref{ODE} are the constant terms for the unlabeled inputs and are the entries on the diagonal of $\boldsymbol{D}_{\mathcal{U}}$. 
Therefore, we can define the vector $\hat{\boldsymbol{b}}$ in \cref{ODE} as having entries $\hat{b}_i={[\boldsymbol{D}_{\mathcal{U}}]}_{i,i}$. 

Finally, we can define \(\boldsymbol{\hat{A}}\) in \cref{ODE} and \(\boldsymbol{F}\) in \cref{constraints}. The matrix \(\boldsymbol{\hat{A}}\) is defined by the equation
		\begin{equation}\label{MatrixA}
			\boldsymbol{\hat{A}}=(\boldsymbol{W}^T-\boldsymbol{D}_\text{out}-\boldsymbol{D}_{\mathcal{V}}){\boldsymbol{X}^\mathcal{T}}^{-1} = (\boldsymbol{W}^T - \boldsymbol{F}_{\text{out}}){\boldsymbol{X}^\mathcal{T}}^{-1}
		\end{equation}
where the diagonal matrix \(\boldsymbol{F}_{\text{out}}\) is defined by the equation $\boldsymbol{F}_{\text{out}} = \boldsymbol{D}_\text{out} + \boldsymbol{D}_{\mathcal{V}}$ and gives the total rate of flux out of each metabolite. 
Similarly, we define the matrix $\boldsymbol{F}_{\text{in}}$ to be the diagonal matrix of all fluxes entering each metabolite, $\boldsymbol{F}_{\text{in}} = \boldsymbol{D}_\text{in} + \boldsymbol{D}_{\mathcal{L}}+\boldsymbol{D}_{\mathcal{U}} $. 
As stated above, the algebraic constraints for constant total metabolite concentrations is that $\boldsymbol{F}_{\text{in}}=\boldsymbol{F}_{\text{out}}$. 
When these constraints are satisfied, we call this diagonal matrix $\boldsymbol{F}$ for simplicity. The diagonal entries of $\boldsymbol{F}$ will be denoted $F_i$ for the flux through the $i^{th}$ metabolite.

\subsection{Scaling for Isotope Enrichment}

Measurements of labeled and unlabeled metabolite concentrations are typically reported as a ``percent enrichment", the percentage of the total pool of a metabolite that is labeled with at least one \(^{13}\)C. 
This is because the mass spectrometer measures a signal intensity at specific masses. 
The intensity of the signal at the mass of the labeled metabolite is divided by the total signal intensity for the metabolite to compute a labeled relative signal intensity. 
Conversion of signal intensity to concentration requires extensive additional experiments to develop calibration curves for each metabolite in the experiment. 
However, the relative signal intensity for metabolite \(i\) is a reasonable approximation of the relative concentration for metabolite \(i\) whether labeled ($x^{\mathcal{L}}_i/x^{\mathcal{T}}_i$) or unlabeled ($x^{\mathcal{U}}_i/x^{\mathcal{T}}_i$). 
Therefore, we scale the model so that the new variables are the proportions of each metabolite that are unlabeled, $\bar{x}_i=x^{\mathcal{U}}_i/x^{\mathcal{T}}_i$ or in matrix form

    \begin{equation}\label{re-scaled variable}
        \boldsymbol{\bar{X}} = {\boldsymbol{X}^\mathcal{T}}^{-1} \boldsymbol{X}^{\mathcal{U}}.
    \end{equation}
 
In this scaled format, it is useful to also write these equations in terms of new parameters. 
Let us define the vector $\vec{\boldsymbol{\alpha}}=\boldsymbol{F}_{\text{in}}^{-1} \hat{\boldsymbol{b}}$ and the matrices $\boldsymbol{B}=\boldsymbol{F}_{\text{in}}^{-1} \boldsymbol{W}^T$,  $\boldsymbol{K}_{\text{in}}={\boldsymbol{X}^{\mathcal{T}}}^{-1} \boldsymbol{F}_{\text{in}}$, and $\boldsymbol{K}_{\text{out}}={\boldsymbol{X}^{\mathcal{T}}}^{-1} \boldsymbol{F}_{\text{out}}$. 
The entries of $\vec{\boldsymbol{\alpha}}$ are $\alpha_i=b_i/{[\boldsymbol{F}_{\text{in}}]}_{i,i}$, 
the proportion of the total flux into metabolite $i$ that comes from an unlabeled input originating outside of the diagram. 
Similarly, for the fluxes originating from other metabolites, the entries of $\boldsymbol{B}$ are $\beta_{i,j}=w_{j,i}/{[\boldsymbol{F}_{\text{in}}]}_{i,i}$, the proportion of the total flux into metabolite $i$ that comes from metabolite $j$. This new notation for the fluxes has been incorporated into a diagram of the network for our example and is shown in \cref{fig:RescaledCyclicPathway}.
The entries of $\vec{\boldsymbol{\alpha}}$ and $\boldsymbol{B}$ are all the dimensionless proportions of fluxes into each metabolite. 
Finally, $\boldsymbol{K}_{\text{in}}$ and $\boldsymbol{K}_{\text{out}}$ are diagonal matrices with entries that can be viewed as turnover rates of the metabolites in units $1/$time. 
The scaled system can be written
    $$
        \boldsymbol{\dot{\bar{X}}} = (\boldsymbol{K}_{\text{in}}\boldsymbol{B}-\boldsymbol{K}_{\text{out}})\boldsymbol{\bar{X}}+\boldsymbol{K}_{\text{in}} \vec{\boldsymbol{\alpha}}.
    $$
As before, if the metabolite concentrations are unchanging then we have the algebraic constraints $\boldsymbol{F}_{\text{in}}=\boldsymbol{F}_{\text{out}}=\boldsymbol{F}$ which gives $\boldsymbol{K}_{\text{in}}=\boldsymbol{K}_{\text{out}}=\boldsymbol{K}$. 
This simplifies the system further to
    \begin{equation}\label{re-scaled ODE}
        \begin{aligned}
            &\boldsymbol{\dot{\bar{X}}} = \boldsymbol{K} (\boldsymbol{B}-\boldsymbol{I})\boldsymbol{\bar{X}}+\boldsymbol{K} \vec{\boldsymbol{\alpha}} \\
            &\boldsymbol{\bar{X}}(0)=\vec{\boldsymbol{1}}
        \end{aligned}
    \end{equation}
where $\vec{\boldsymbol{1}}$ is the vector with all entries equal to one. 

This formulation of the problem is powerful because it allows us to write the model variables entirely in terms of proportions unlabeled. 
This proportion can be measured directly as one minus the isotopic enrichment. 
Additionally, the parameters in $\vec{\boldsymbol{\alpha}}$ and $\boldsymbol{B}$ are also ratiometric. 
They are proportions of each flux to total flux through the metabolites so we will refer to them as flux ratio parameters. 
The other parameters in $\boldsymbol{K}$ are the turnover rates for each metabolite. 
There is no explicit dependence on the total concentrations. 
This means we can construct the model and fit the parameters in this format without measuring metabolite concentrations. 
To convert the resulting values for the parameters back to fluxes in units of concentration per unit of time would still require the costly measurement of total concentrations. 
However, in many instances functional changes in metabolic networks could be gleaned from changes in the turnover rates and relative flux values which are the parameters in equation \cref{re-scaled ODE}. 

\subsection{Number of Independent Parameters}\label{NumberParameters}
 
Since the computational time for parameter estimation can quickly become impractical, it is essential to reduce the number of parameters as much as possible. 
In this section, we use the algebraic constraints to find the number of independent parameters which must be fit. 
We begin with a general statement of the number of parameters needed. 
Next, we investigate the special case where each metabolite has a non-zero flux exiting the graph. 
Finally, we discuss the more difficult case, where some exiting fluxes are missing from the reaction diagram.

\begin{theorem}
The algebraic constraints in $\cref{constraints}$ can be used to rewrite the system in \cref{ODE} in terms of $R-N$ parameters.
\end{theorem}

\begin{proof}

In the original formulation before scaling, the model contains parameters for the fluxes $\vec{\boldsymbol{f}}={f_1, f_2, \cdots f_R}$. 
These fluxes can be organized into an $N \times R$ incidence matrix, $\boldsymbol{M}$. 
In the context of chemical reactions, this is equivalent to the stoichiometry matrix for this set of reactions. 
The rows are the metabolites and the columns correspond to the fluxes. 
The entries of $\boldsymbol{M}$, $m_{i,j}$ will be $1$ if flux $j$ is an input for metabolite $i$, and $-1$ if it is an output from metabolite $i$. 
All other entries will be zero. 
Our algebraic constraints that the total input flux must equal the total output flux for each metabolite can therefore be written
    $$
        \boldsymbol{M} \vec{\boldsymbol{f}}=0
    $$
The matrix $\boldsymbol{M}$ will be full rank. 
This is easiest to observe in the row reduction of $\boldsymbol{M}^T$. 
We can assume that there is at least one input since we must have labeled input. 
The row of $\boldsymbol{M}^T$ for this input will include a $1$ in the column for the metabolite entered, and a zero for all other columns. 
We can therefore eliminate all other entries in this column. 
Any flux that connected this metabolite to another metabolite will now have a row with only one non-zero entry. 
These can be used to eliminate other entries in those columns, and so on. 
Because the graph is connected this will eventually reach all columns leaving a pivot in each column. 
Since the matrix $\boldsymbol{M}$ has rank $N$, we can solve for $R-N$ of the fluxes in terms of the other $N$. 
In general, the system will have $R-N$ independent parameters.
\end{proof}

\begin{figure}[H]
        \centering
        \includegraphics[width=0.8\linewidth]{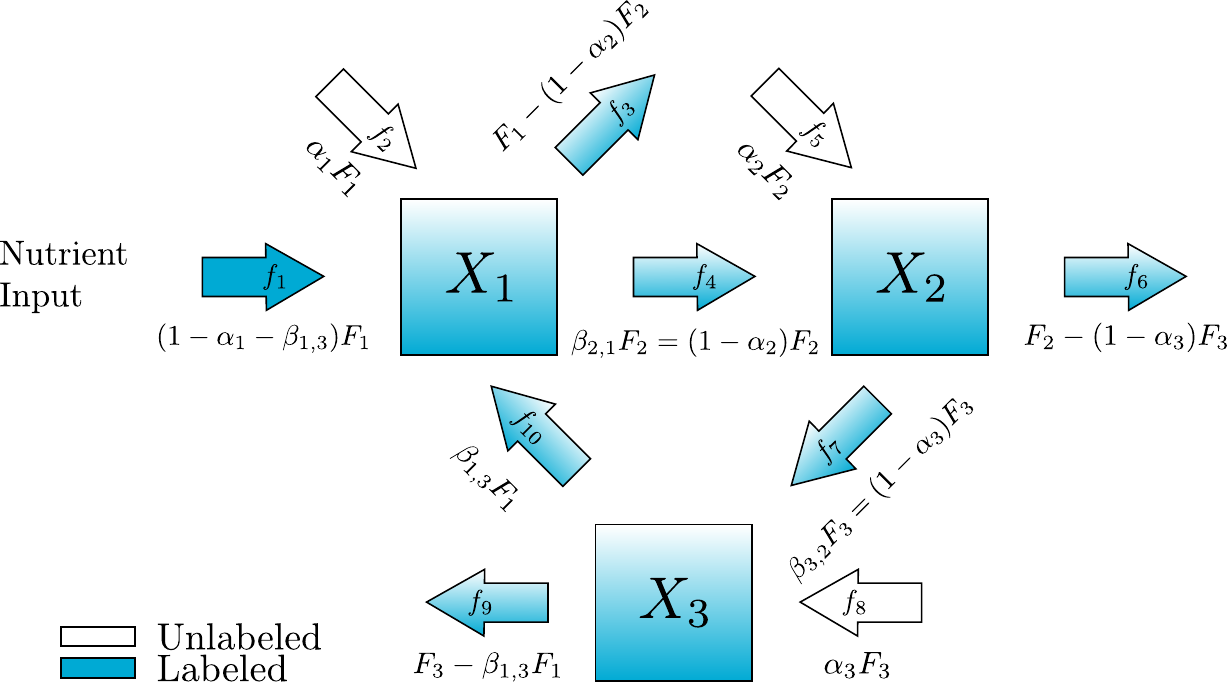}
        \caption{An example demonstrating the new notation with proportional rate parameters. This is the same three metabolite pathway shown in \cref{Cyclic} }
        \label{fig:RescaledCyclicPathway}
    \end{figure}

Before we turn to the special case, let us begin by counting the parameters in our scaled system \cref{re-scaled ODE}. 
This model has $N$ parameters corresponding to the non-zero entries of $\boldsymbol{K}$, $k_i$ for the turnover rate of each metabolite, $|\mathcal{E}_\mathcal{U}|$ parameters for the non-zero entries of $\vec{\boldsymbol{\alpha}}$ for the unlabeled inputs, and $|\mathcal{E}_\mathcal{W}|$ parameters corresponding to the non-zero entries of $\boldsymbol{B}$ for the fluxes between metabolites. 
Since $\alpha_i$ and $\beta_{i,j}$ are all of the dimensionless proportions of fluxes into metabolite $i$, for any metabolite which has no labeled inputs, we know that
    \begin{equation}\label{proportions in equal 1}
        \sum_{j=1}^N \beta_{i,j} + \alpha_i =1.
    \end{equation}
We can therefore, without loss of generality, replace the first non-zero entry in these rows 
    \begin{equation}\label{beta substitution}
        \beta_{i,k}=1-\alpha_i -\sum_{j=k+1}^N \beta_{i,j}.
    \end{equation}
There are therefore $|\mathcal{E}_{\mathcal{W}}|-(N-|\mathcal{E}_{\mathcal{L}}|)$ remaining independent $\beta$ values. 
In our example, this means we can replace $\beta_{2,1}$ with $(1-\alpha_2)$ and $\beta_{3,2}$ with $(1-\alpha_3)$ so we have written both forms in \cref{fig:RescaledCyclicPathway}. Now we can write our system in terms of the seven remaining parameters. 

    \begin{equation*} 
        \boldsymbol{\dot{\bar{X}}}= \begin{pmatrix}
        k_1 & 0 & 0 \\ 0 & k_2 & 0 \\
        0 & 0 & k_3
        \end{pmatrix}
        \left(\begin{pmatrix}
        -1 & 0 & \beta_{1,3}\\
        (1-\alpha_2) & -1 & 0\\
        0 & (1-\alpha_3) & -1
        \end{pmatrix} \boldsymbol{\bar{X}} +
        \begin{pmatrix}
        \alpha_1 \\ \alpha_2\\ \alpha_3
        \end{pmatrix}  
        \right)
    \end{equation*}

A graph which is an arborescence will have $|\mathcal{E}_{\mathcal{L}}|=1$ and $|\mathcal{E}_{\mathcal{W}}|=N-1$. Therefore, such a graph will have no remaining independent $\beta$ parameters. On the other hand, a graph that is not an arborescence will have either $|\mathcal{E}_{\mathcal{L}}|>1$ or $|\mathcal{E}_{\mathcal{W}}|>N-1$ or both. Therefore, any graph which is not an arborescence will have non-zero $\beta$ parameters in this scaling. 

So far we have counted the turnover rate parameters, $k_i$, and the proportional flux parameters, $\alpha_i$ and $\beta_{i,j}$. This brings us to a total number of parameters in \cref{re-scaled ODE} of $N+|\mathcal{E}_\mathcal{U}|+|\mathcal{E}_{\mathcal{W}}|-(N-|\mathcal{E}_{\mathcal{L}}|)=R-|\mathcal{E}_{\mathcal{V}}|$. 
In developing this parameter set, we have labeled the edges in $\mathcal{E}_{\mathcal{U}}$ with $\alpha$s and those in $\mathcal{E}_{\mathcal{W}}$ with $\beta$s. 
Let us now turn our attention to $\mathcal{E}_{\mathcal{L}}$ and $\mathcal{E}_{\mathcal{V}}$. 
Because of our choice to scale all inputs to the total input for each metabolite, a labeled input edge in $\mathcal{E}_{\mathcal{L}}$ can also scale to the influx into its node $i$ as $(1-\alpha_i -\Sigma_{j=1}^N \beta_{i,j})F_{i}$. 
In our example, this means that the labeled input will be $f_1=(1-\alpha_1-\beta_{1,3})F_{1}$. 
In practice, it is typical to have only one labeled input so each of our examples has one edge in $\mathcal{E}_{\mathcal{L}}$ which is an input into the first metabolite. 

The labeling of the edges in $\mathcal{E}_{\mathcal{V}}$ may pose a new complication because of our algebraic constraints. 
Notice first that in the special case $|\mathcal{E}_{\mathcal{V}}|=N$, the number of parameters already defined after the substitution in \cref{beta substitution} will be exactly the minimum number, $R-N$. 
This is because in this case, every node has an output vector that can be used to satisfy the algebraic constraint that the total output from each node must equal the flux through that node. 
Since the total of all fluxes out of node $i$ must be $F_{i}$, the flux from metabolite $i$ exiting the diagram will be $F_{i}-\sum_{j \ne i} \beta_{j,i} F_{j}$. 
In our example, this means that $f_3=F_{1}-(1-\alpha_2)F_{2}$, $f_6=F_{2}-(1-\alpha_3) F_{3}$, and $f_{10}=F_{3}-\beta_{1,3}F_{1}$. 
Since every node in the example has an output, we have completed our parameterization in a way that gives the minimum number of parameters and does not depend on any concentrations which would require additional measurements.  

We will finish this section with a discussion about the case where $\mathcal{E}_{\mathcal{V}} < N$. 
In this case, there is at least one metabolite that does not have an outflux leaving the pathway. 
This means that we still have too many parameters in our formulation. 
However, our algebraic constraints now imply that the equation for the missing outflux must be zero. 
That is, for any node $i$ which does not have a flux out of the pathway, $F_{i}=\Sigma_{j \ne i} \beta_{j,i} F_{j}$. 
This constraint is challenging because it is a relationship among the total fluxes through multiple metabolites. 
It can be solved for one of the $\beta$ parameters or one of the $k$ parameters, however, the resulting equation will depend on the total concentrations of node $i$ and any immediately downstream nodes. 
In our example, suppose there is no outflux from metabolite $3$. 
Then $f_{9}=0$ implies $F_{3}=\beta_{1,3}F_{1}$ To write this in our preferred parameters we have $\beta_{1,3}=\frac{k_{3}}{k_{1}} \frac{x_3^{\mathcal{T}}}{x_1^{\mathcal{T}}}$ or $k_3=\beta_{1,3} k_1 \frac{x_1^{\mathcal{T}}}{x_3^{\mathcal{T}}}$. 
In general, any metabolite that does not have an outflow will create a constraint on the parameters which requires reintroduction of the total concentrations of the metabolites. 
In designing KFP experiments, this leads to a decision that must be made. 
In some cases, the experimenter may already plan to conduct additional experiments to measure the concentrations of the relevant metabolites, either as part of their larger study to understand mechanistic changes or because they wish to convert the relative fluxes and turnover rates in equation \cref{re-scaled ODE} back to the original fluxes in \cref{ODE}.
In this case, the constraint induced by the missing outflux can be used to further reduce the number of parameters. If the total concentrations will not be measured, we will need to include a flux out of the pathway for every metabolite. 
The additional parameter not only increases computation time, but it may also increase the number of time points needed to accurately and confidently estimate relative fluxes and turnover rates. 

In this section, we have used the structure of the directed weighted graph representation of a metabolic pathway to formulate differential equations for the unlabeled metabolite concentrations. 
Scaling this model in terms of isotope enrichment allows us to write the model entirely in terms of relative fluxes and turnover rates. 
Using the algebraic constraints that the flux in and out of each metabolite must be equal, we can further reduce the number of parameters. 
In the case where every metabolite has a flux out of the pathway ($|\mathcal{E}_{\mathcal{V}}|=N$) the model can be written in terms of the minimal number of parameters using only relative fluxes and turnover rates. 
In the case where not all metabolites have a flux out of the pathway, the resulting minimal set of parameters will include dependence on the total concentrations. 

\section{Steady State Problem}\label{sec3}

In this section, we turn our attention to the second question. 
What can be learned from the steady-state fractional isotopic enrichments? Other methods such as \(^{13}\)C metabolic flux analysis use data collected only after the isotopic steady state is reached. 
In practice, one must design the experimental protocol without preexisting knowledge of the rates of the reaction. 
At the same time, collecting data from time points very soon after the isotope switch without perturbing cellular metabolism can be technically challenging. 
Therefore, the metabolites may have attained their isotopic steady state by the time the first data points are collected. 
We show below that only the proportional flux parameters can be estimated from the steady-state information and only under specific conditions on the graph structure. 

\subsection{Proportional Flux Parameters and the Isotopic Steady State}
\begin{theorem}
    The steady state of \cref{re-scaled ODE} is unique. This steady state depends on the proportional flux parameters and does not depend on the turnover rates.
\end{theorem}

\begin{proof}
The steady state of \cref{re-scaled ODE}, \(\boldsymbol{\bar{X}}^{ss}\), is found by solving the matrix equation

    \begin{equation}\label{re-scaled ODE SS}
        (\boldsymbol{I}-\boldsymbol{B})\boldsymbol{{\bar{X}}^{ss} = \vec{\boldsymbol{\alpha}}}
    \end{equation}
The matrix \(\boldsymbol{K}\) in equation \cref{re-scaled ODE} does not appear in \cref{re-scaled ODE SS} because $\boldsymbol{K}$ is invertible. 
By definition, the diagonal entries of $\boldsymbol{K}$, \(F_i/ X^{\mathcal{T}}_i > 0 \) are non-zero because the metabolite must be present (\(X^{\mathcal{T}}_{i} > 0\)) and there must be flux through it (\(F_i > 0\)). 

Next we will show that $(\boldsymbol{I}-\boldsymbol{B})$ is weakly chained diagonally dominant. Because the entries of the matrix $\boldsymbol{B}$ in row $i$ are the proportions of the fluxes into metabolite $i$, the row sums of $\boldsymbol{B}$ are between zero and one. Subtracting, this implies that $(\boldsymbol{I}-\boldsymbol{B})$ has nonnegative row sums and is therefore weakly diagonally dominant. Now, we must use two properties of the underlying graph structure, first that there is at least one labeled input, and second that there is a path to every node from a labeled input. The row of $(\boldsymbol{I}-\boldsymbol{B})$ corresponding to a node that receives a labeled input will be strictly diagonally dominant. Since there is a path to all other nodes in the pathway, a directed graph constructed from $(\boldsymbol{I}-\boldsymbol{B})$ will have a path back to the labeled input node.  Since  $(\boldsymbol{I}-\boldsymbol{B})$ is weakly chained diagonally dominant, it is invertible \citep{shivakumar1974sufficient}.
Therefore, the steady-state values \(\boldsymbol{\bar{X}}^{ss}=(\boldsymbol{I}-\boldsymbol{B})^{-1} \vec{\boldsymbol{\alpha}}\) only depends on the proportional flux parameters, \(\alpha_i\) and \(\beta_{i,j}\).
\end{proof}

Note that in the special case where each metabolite receives at least one input from $\mathcal{E}_{\mathcal{L}}$ or $\mathcal{E}_{\mathcal{U}}$ we can show that $(\boldsymbol{I}-\boldsymbol{B})$ has strictly positive row sums. This is equivalent to the condition that $\alpha_i>0$ for all metabolites with no labeled input. Consider a metabolite $i$ with no labeled input. Recall from \cref{proportions in equal 1} that the entries of row $i$ of matrix $\boldsymbol{B}$ must add to $1-\alpha_i$. This implies that the row sum for the $i^{th}$ row of $(\boldsymbol{I}-\boldsymbol{B})$ is $\alpha_i$. In the case of a row corresponding to a metabolite which does have labeled input, the row sum of $\boldsymbol{B}$ must be strictly less than one since we require a non-zero labeled influx. Therefore, $(\boldsymbol{I}-\boldsymbol{B})$ is a Z-Matrix with positive row sums and is invertible \citep{berman1994nonnegative}.

We have chosen the scaled parameters so that the equilibrium equation \cref{re-scaled ODE SS} will have no dependence on the turnover rates in $\boldsymbol{K}$, therefore these turnover rates can not be estimated without some data collected during the transient time. While we have defined the function from $\alpha$ and $\beta$ parameters to the steady state values, this function is not, in general, invertible. Recall from \cref{NumberParameters} that the number of nonzero parameters in $\vec{\boldsymbol{\alpha}}$ is $|\mathcal{E}_{\mathcal{U}}|$ and that the number of non-zero parameters in $\boldsymbol{B}$ is $|\mathcal{E}_{\mathcal{W}}|$. After substitution using \cref{beta substitution}, there will be $|\mathcal{E}_{\mathcal{L}}|+|\mathcal{E}_{\mathcal{U}}|+|\mathcal{E}_{\mathcal{W}}|-N$ proportional flux parameters in \cref{re-scaled ODE SS}. 
The simple condition that there must be at least as many steady state values as parameters leads to the following necessary condition for recovery of the $\alpha$ and $\beta$ parameters 

\begin{equation}\label{Necessary Condition SS}
    |\mathcal{E}_{\mathcal{L}}|+|\mathcal{E}_{\mathcal{U}}|+|\mathcal{E}_{\mathcal{W}}| \le 2N.
\end{equation}

Note that the number of arrows leaving the pathway, $\mathcal{E}_{\mathcal{V}}$ is not included in this condition. As we saw in \cref{NumberParameters}, if all nodes have an output vector in $\mathcal{E}_{\mathcal{V}}$ then these output vectors are determined by the algebraic constraints. If one of the nodes does not have such an output vector, the algebraic constraint gives a new condition on parameters. However, this can not reduce the number of parameters in $\vec{\boldsymbol{\alpha}}$ or $\boldsymbol{B}$ so can be ignored for this discussion.

A few simple examples illustrate why this condition \cref{Necessary Condition SS} is not sufficient. If $\vec{\boldsymbol{\alpha}}=\vec{\boldsymbol{0}}$ then the only steady state is all zeros. Intuitively, this case has no unlabeled input so at steady state, all metabolites are fully labeled. Clearly no information about the values of the $\beta$ parameters can be deduced from this steady state. Another example is if metabolite $i$ only receives input from metabolite $j$, that is $\alpha_i=0$, $\beta_{i,k}=0$ for $k \ne j$. In this case, the metabolites $i$ and $j$ will approach the same steady state and the steady states will provide fewer than $N$ equations for specifying underlying flux ratios. 

In the case where the graph is an arborescence, the condition is met and the $\alpha$ proportional flux parameters can always be recovered from steady state values. In this case the steady state of the root node will be its proportion of unlabeled input. Every other node $i$, can be written in terms of only its unique upstream metabolite $j$. $\bar{x}_i^{ss}=1-(1-\alpha_i)(1-\bar{x}^{ss}_j$). Therefore $\alpha_i=\frac{\bar{x}_i^{ss}-\bar{x}_j^{ss}}{1-\bar{x}_j^{ss}}$.

Next we will look in detail at two simple examples which illustrate this necessary condition. In the first, there are two metabolites and the reaction between them is not reversible. In this case, the condition \cref{Necessary Condition SS} is met and the parameters can be written as functions of the steady state values. Using Bayesian parameter estimation of simulated noisy data, we see that the $\alpha$ parameters are recovered accurately and confidently. Next, we consider the similar case with a reversible reaction. In this case the \cref{Necessary Condition SS} condition is not met and the estimation of parameters is problematic. 

\subsection{Examples: Irreversible and Reversible Reactions}

Given data collected only at the steady state of the system, proportional flux parameters can only be uniquely recovered from steady-state estimates if the map from proportional flux parameter values to steady-state values is one-to-one. In this section we explore two examples, one irreversible and one reversible which demonstrate this condition.

\begin{figure}[H]
	\centering
	\includegraphics[width=0.70\textwidth]{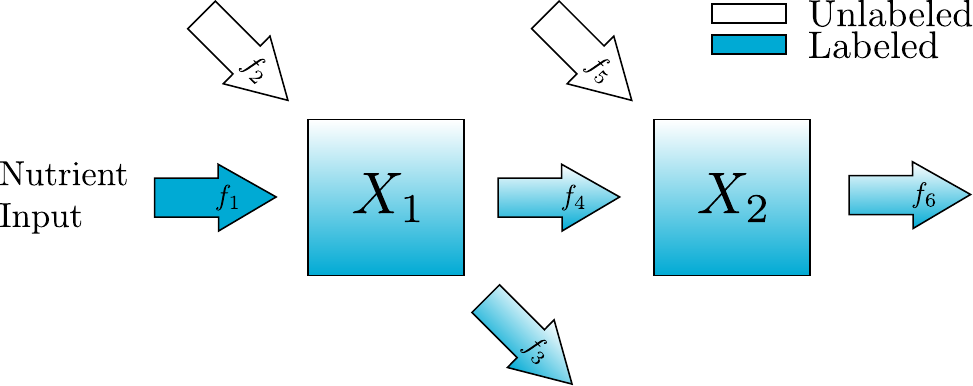}
	\caption{Irreversible two-metabolite model in metabolic steady state. Isotopic label feeds into metabolite \(X_1\) through the nutrient input, $f_1$, while unlabeled metabolite enters through an external source, $f_2$. A mix of labeled and unlabeled metabolite is converted into metabolite \(X_2\), $f_4$, or leaves the scope of the diagram, $f_3$. Unlabeled \(X_2\) enters the system through another external source, $f_5$, while a mix of labeled and unlabeled \(X_2\) leaves the scope of the diagram, $f_6$}
    \label{fig:One direction 2 metabolites}
\end{figure}

The scaled equation for \cref{fig:One direction 2 metabolites} is
    \begin{equation}\label{re-scaled 2M ODE}
        \boldsymbol{\dot{\bar{X}}} =\begin{pmatrix}k_1 & 0 \\ 0 & k_2\end{pmatrix} \left(\begin{pmatrix}-1 & 0 \\ (1-\alpha_2) & -1 \end{pmatrix}\boldsymbol{\bar{X}} + \begin{pmatrix}\alpha_1 \\ \alpha_2\end{pmatrix}\right)
    \end{equation}
which has a steady state 
    \begin{equation}\label{irreversible SS}
        \boldsymbol{\bar{X}}^{ss} = \begin{pmatrix}
        1 & 0 \\ (1-\alpha_2) & 1\end{pmatrix}\begin{pmatrix}
            \alpha_1 \\ \alpha_2
        \end{pmatrix} = \begin{pmatrix}
            \alpha_1 \\ \alpha_1 - \alpha_1\alpha_2 + \alpha_2
        \end{pmatrix}
    \end{equation}

Suppose we can get estimates for the steady-state values of the metabolites, \((\bar{x}_1^{ss}, \bar{x}_2^{ss})\). The proportional parameters \(\alpha_1\) and \(\alpha_2\) can be recovered with the map 
    \[(\bar{x}_1^{ss},\bar{x}_2^{ss}) \to \left(\bar{x}_1^{ss}, \frac{\bar{x}_2^{ss}-\bar{x}_1^{ss}}{1-\bar{x}_1^{ss}}\right) = (\alpha_1, \alpha_2)\]

To illustrate the ability to recover $\alpha_1$ and $\alpha_2$ in this example, we have used Bayesian parameter estimation with noisy simulated data as described in Appendix \cref{appendix: Data and BPE}. In all of our Bayesian parameter estimation examples we have used the following general procedure. We generate solutions to the ODE system using known parameter values. Then at each time point, we produce three noisy data points to simulate technical replicates of the experiment. To explore the role of the noise, we have added normally distributed noise to the true solution with a standard deviation that is 2.5, 5, or 10\% of the true value. To explore the dependence on the number of time points we have included 3, 5, or 10 equally spaced time points in each parameter estimation. All parameters, including proportional flux parameters and turnover rates are given naive priors. In this section we show only the posterior distributions for the proportional flux parameters, since we anticipate poor estimates of turnover rates. The simulated data, the estimates for the turnover rates, and the credible intervals for each parameter estimation can all be found in the supplement. 

In the violin plots in \cref{fig:2MSimilarDenAB}, we see that the method consistently provides accurate recovery of both $\alpha_1$ and $\alpha_2$ as the mode of the posterior distributions. The spread of the violin plots or box plots, is a visual indicator of the uncertainty in our estimate. The comparison across number of time points shows, not surprisingly, that more time points generally provides a more accurate and less variable estimate. The results are relatively insensitive to the level of noise in the data. 
Increasing the noise level only slightly increases the error and the uncertainty in the estimate. 
 
    \begin{figure}[H]
        \centering
        \includegraphics[width=0.9\textwidth]{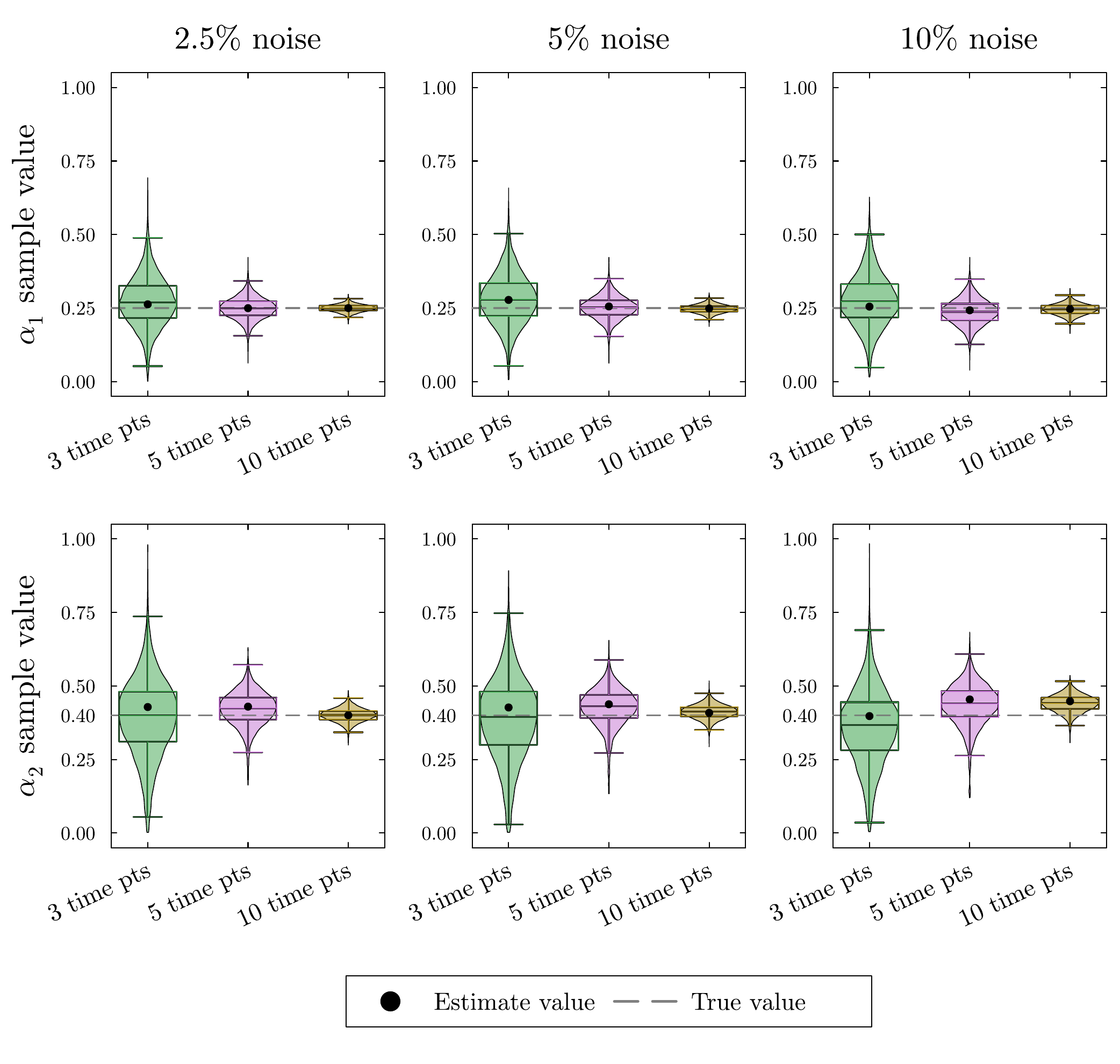}
        \caption{Violin and box plots for posterior distributions of \(\alpha_1\) and \(\alpha_2\) for the irreversible two-metabolite model with similar turnover rates. The simulated data are taken from solutions of system \cref{re-scaled 2M ODE} with parameters \(k_1 = 7/20\), \(\alpha_1 = 1/4\), \(k_2 = 3/10\), \(\alpha_2 =2/5\) at 3, 5 or 10 equally distributed time points. At each time point, three noisy data points are simulated by adding random normally distributed noise to the true solutions with a standard deviation of 2.5, 5, or 10\% of the true value. Since \(\alpha_1\), and \(\alpha_2\) are proportions, we use the natural naive prior distributions
		$\alpha_1 \sim U(0,1),\; \alpha_2 \sim U(0,1)$. The mode of the posterior samples is shown as a dot and is taken as the estimated value. Outliers are excluded in the box plots}
        \label{fig:2MSimilarDenAB}
    \end{figure}

\begin{figure}[H]
	\centering
	\includegraphics[width=0.70\textwidth]{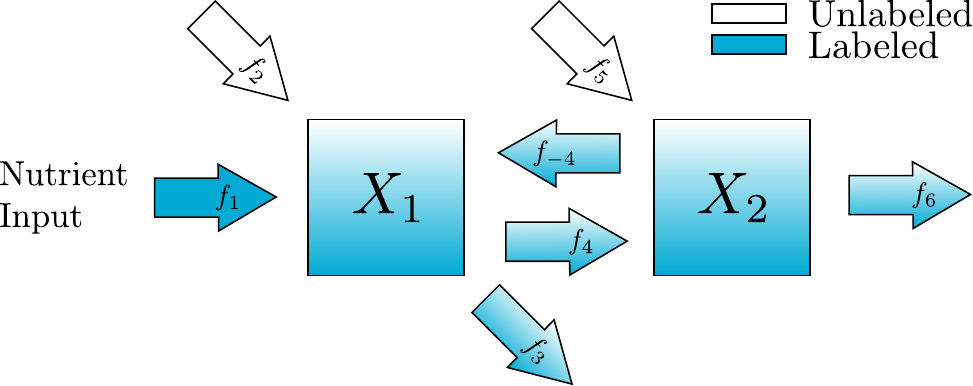}
	\caption{Reversible two-metabolite model in metabolic steady state with reversible reaction. Isotopic label feeds into metabolite \(X_1\) through the nutrient input, $f_1$, while unlabeled metabolite enters through an external source, $f_2$. A mix of labeled and unlabeled metabolite is converted into metabolite \(X_2\), $f_4$, or leaves the scope of the diagram, $f_3$. Unlabeled \(X_2\) enters the system through an external source, $f_5$, while a mix of labeled and unlabeled \(X_2\) can return to metabolite \(X_1\), $f_{-4}$, or leave the scope of the diagram, $f_6$}
    \label{fig:Reversible 2 metabolites}
\end{figure}

Next, we turn to the second example shown in \cref{fig:Reversible 2 metabolites}. In this reversible reaction, the number of proportional flux parameters is 3, and therefore it does not satisfy the condition \cref{Necessary Condition SS}. The scaled system for this example is
    \begin{equation}\label{re-scaled R2M ODE}
        \boldsymbol{\dot{\bar{X}}} = \begin{pmatrix}k_1 & 0 \\ 0 & k_2\end{pmatrix}\left( \begin{pmatrix}-1 & \beta_{1,2} \\ (1-\alpha_2) & -1 \end{pmatrix}\boldsymbol{\bar{X}} + \begin{pmatrix}\alpha_1 \\ \alpha_2\end{pmatrix}\right)
    \end{equation}
which has a steady state 
    \begin{equation}\label{reversible SS}
        \boldsymbol{\bar{X}}^{ss} = \frac{1}{1-\beta_{1,2}(1-\alpha_2)}\begin{pmatrix}
        1 & \beta_{1,2} \\ (1-\alpha_2) & 1\end{pmatrix}\begin{pmatrix}
            \alpha_1 \\ \alpha_2
        \end{pmatrix} = \begin{pmatrix}
            \frac{\alpha_1+\beta_{1,2}\alpha_2}{1-\beta_{1,2}(1-\alpha_2)} \\ \frac{\alpha_1 - \alpha_1\alpha_2 + \alpha_2}{1-\beta_{1,2}(1-\alpha_2)}
        \end{pmatrix}
    \end{equation}

Suppose we can get estimates for the steady-state values of the metabolites, \((\bar{x}_1^{ss}, \bar{x}_2^{ss})\). Since these values depend on all three parameters, there is no way to recover unique values for the parameter with only two pieces of information.

\begin{figure}[H]
        \centering
        \includegraphics[width=0.85\textwidth]{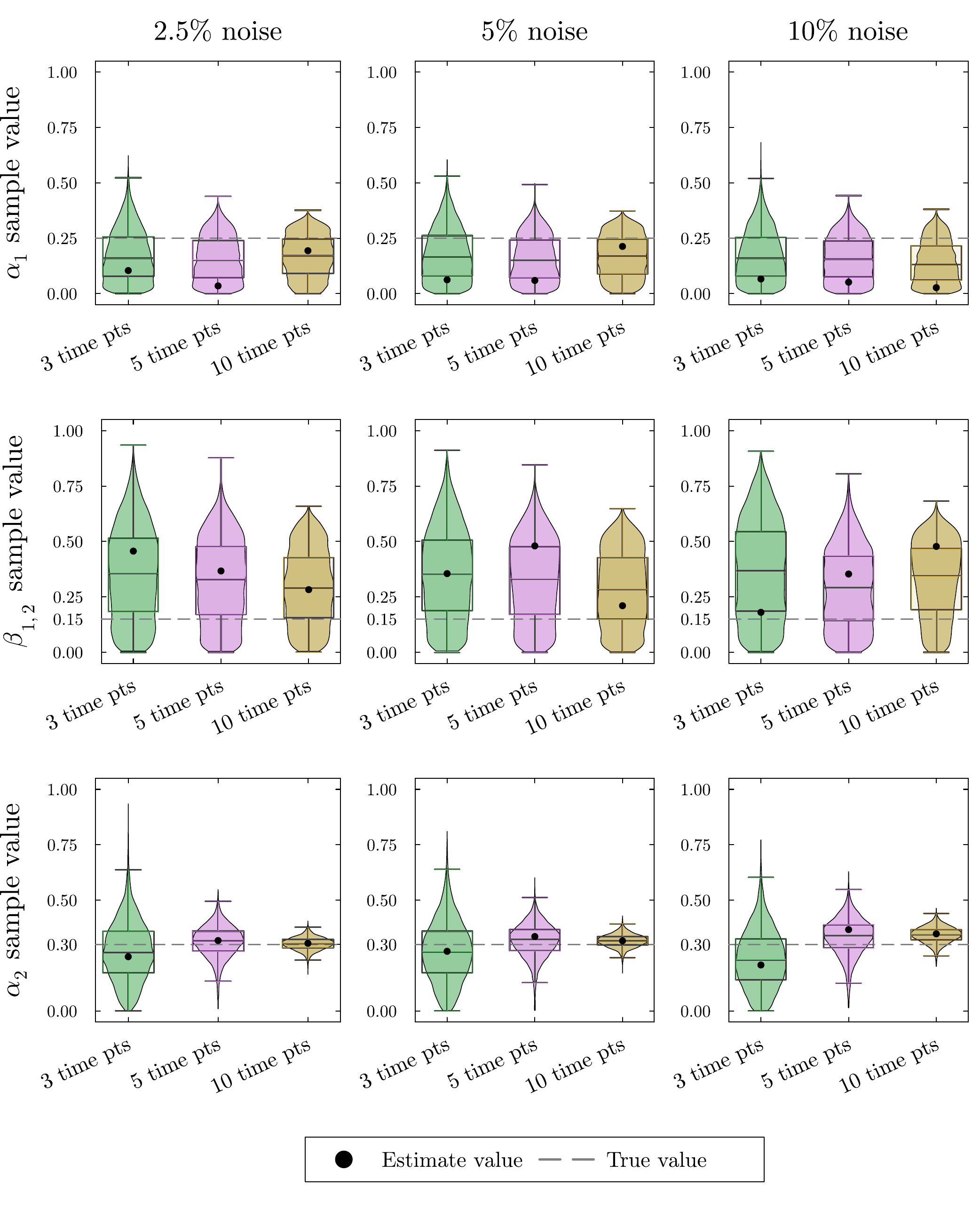}
        \caption{Violin and box plots for posterior distributions of \(\alpha_1\), \(\beta_{1,2}\), and \(\alpha_2\) for the reversible two-metabolite model with similar turnover rates. The simulated data are taken from solutions of system \cref{re-scaled R2M ODE} with parameters \(k_1 = 7/20\), \(\alpha_1 = 3/10\), \(\beta_{1,2} = 3/20\) \(k_2 = 3/10\), \(\alpha_2 =1/4\) at 3, 5 or 10 equally distributed time points. At each time point, three noisy data points are simulated by adding random normally distributed noise to the true solutions with a standard deviation of 2.5, 5, or 10\% of the true value. Since \(\alpha_1\), and \(\alpha_2\) are proportions, we use the natural naive prior distributions
		$\alpha_1 \sim U(0,1),\;\beta_{2,1} \sim U(0,1), \; \alpha_2 \sim U(0,1)$. The mode of the posterior samples is shown as a dot and is taken as the estimated value. Outliers are excluded in the box plots}
        \label{fig:R2MSimilarDenABG}
    \end{figure}

Providing Bayesian parameter estimation with simulated data of this system as described in Appendix A shows the inability to recover all parameter values accurately.
 In the violin plots in \cref{fig:R2MSimilarDenABG}, we see that we only get moderately accurate estimates for \(\alpha_2\). 
As before, as we increase the number of time points, the uncertainty for \(\alpha_2\) estimates decreases and the results are relatively insensitive to the magnitude of the noise. 
The estimates for \(\alpha_1\) and \(\beta_{1,2}\) on the other hand are completely unsatisfactory. The estimates of \(\alpha_1\) are much lower than the true value, while the estimates for \(\beta_{1,2}\) are much higher than the true value. 
The sample distributions for both \(\alpha_1\) and \(\beta_{1,2}\) show that there is a high uncertainty in the estimates. 
Since the two steady-state values depend on the values of \(\alpha_1\), \(\beta_{1,2}\), and \(\alpha_2\), several combinations of parameter values give the same steady-state value. While the method may be able to correctly determine combinations of these parameters that give the correct steady-state values for both \(\bar{x}_1\) and \(\bar{x}_2\), it cannot discern the correct parameter values.  

We further hypothesize that $\alpha_1$ and $\beta_{1,2}$ are more difficult to estimate because of the structure of our steady state solutions. The true values for  \(\alpha_1\), \(\beta_{1,2}\), and \(\alpha_2\) are all fairly small, so if we ignore the quadratic terms in \cref{reversible SS} we obtain approximate steady states $\left(\frac{\alpha_1}{1-\beta_{1,2}} , \frac{\alpha_1 + \alpha_2}{1-\beta_{1,2}}\right)$. In this expression, $\alpha_1$ and $\beta_{1,2}$ always occur in the same combination and therefore cannot be distinguished. 
 
In this section we have shown that the estimation of parameters may be severely limited if only steady state data is obtained. In this case, the values for the turnover rates can not be estimated as they do not appear in the steady state equation. It may, however, be possible to estimate the proportional flux parameters if the number of measured steady state exceeds the number of parameters. The specific example with irreversible and reversible reactions illustrates this necessary condition. The examples also makes clear the utility of using a Bayesian parameter estimation. A classic least squares parameter estimate will be able to replicate the limited data but will not necessarily give accurate estimates of the underlying fluxes. The Bayesian approach, however, allows the researcher to simultaneously quantify the best estimate and the certainty in that estimate. 
This allows for iterative improvement of the experimental protocol for a particular pathway and lowers the chances of misinterpretation.

\section{Fast-Slow Analysis}\label{sec4}

In this section, we address separation in time scales for different metabolites. Here, we will be comparing turnover rates which are the ratio of the total flux through a metabolite to the total concentration of that metabolite. A metabolite may be considered fast if the fluxes into and out of that metabolite are large, or if the total concentration of that metabolite is small. Here we will begin with the general proof that data collected from a reaction network with a disparity of time scales does not contain fast dynamics if the metabolite with the fast turnover rate is not receiving pure isotopically labeled nutrient. Next we will investigate a chain reaction network with the fast metabolite either first or later in the chain. Finally, we will return to the reversible two metabolite model and use Bayesian parameter estimation to illustrate the accuracy and reliability of estimates of the turnover rates. 

\subsection{Fast-Slow Dynamics Limit Turnover Rate Identifiability}\label{sec4.1}

Here we show that if the metabolite with the fast turnover rate is not the metabolite receiving the labeled nutrient input, the solution is entirely on the slow manifold and the parameter for the fast turnover rate cannot be estimated.

    \begin{theorem}\label{FastSlowTheorem}
        In the scaled reaction network model in \cref{re-scaled ODE}, suppose \(k_n \gg k_i\) \(\forall i \neq n\). If \(X_n\) is not the target of any edge in \(\mathcal{E}_\mathcal{L}\), then the initial condition is on the slow manifold.
    \end{theorem}

    \begin{proof}\label{Thm1proof}
        Let \(k_n\gg k_i\) \(\forall i\neq n\). Then, the slow manifold is obtained by setting $\dot{\bar{x_n}}=0$ 
            \begin{equation}\label{Thm Slow Manifold}
                \sum_{i=1}^N(\boldsymbol{B}-\boldsymbol{I})_{n,i}\bar{x}_i + \alpha_n = 0.
            \end{equation}
        Since \(X_n\) is not the target of any edge in \(\mathcal{E}_\mathcal{L}\), the sum of the proportional flux parameters coming into this metabolite will be one.  By \cref{proportions in equal 1} we have 
            \[\sum_{i=1}^{N}\beta_{n,i}+\alpha_n = 1.\]
       Plugging the initial condition,  \(\bar{x}_i = \vec{1}\) into the equation for the slow manifold, we have 
            \[\sum_{i=1}^N(\boldsymbol{B}-\boldsymbol{I})_{n, i} + \alpha_n = \sum_{i=1}^{N}\beta_{n,i} - 1 + \alpha_n = 0.\]
        Hence, the initial condition is on the slow manifold in \cref{Thm Slow Manifold}.
    \end{proof}

As a result of \cref{FastSlowTheorem}, fast dynamics will not be present in the system and the experimenter will not be able to estimate the faster turnover rate \(k_n\). 
In this situation, the experimenter may still be able to estimate all the slower turnover rates and the proportional flux paramters.

In the special case where the metabolite with the quickest turnover rate \(X_n\) is a target of an edge in \(\mathcal{E}_\mathcal{L}\), the fast subsystem is defined by
    \begin{equation*}\label{Genereal fast system}
        \dot{\bar{x}}_n = k_n\left(-\bar{x}_n+\alpha_n+\sum_{\substack{i=1 \\ i\neq n}}^N \beta_{n,i}\right).
    \end{equation*}
Fast dynamics are present in this system but estimating \(k_n\) relies on the ability of the experimenter to collect time points on the scale of \(k_n\).

In the following example, we show that data collected from a chain reaction network with a disparity of time scales does not contain fast dynamics if the slow turnover rate is upstream of a more rapid turnover rate. 
We will begin by formulating the general reversible chain reaction network model and its fast and slow subsystems. 

\subsection{Chain Reaction Model}\label{sec4.2}

\begin{figure}[H]
	   \centering
	   \includegraphics[width=0.85\textwidth]{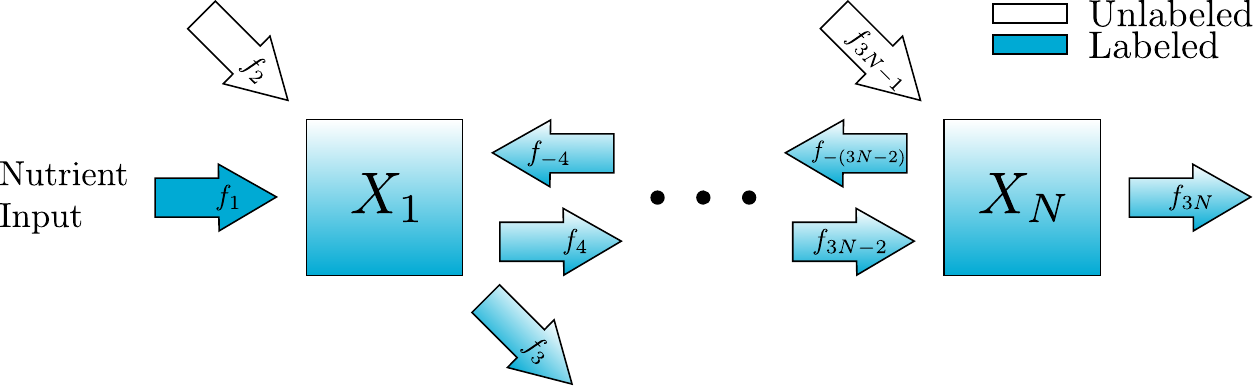}
	   \caption{\(N\)-metabolite reversible chain reaction pathway in metabolic steady state. Isotopic label feeds into metabolite \(X_1\) through the nutrient input while unlabeled metabolite enters through an external source. A mix of labeled and unlabeled metabolite is converted into the following metabolite, \(X_2\), or leaves the scope of the diagram. Unlabeled \(X_2\) enters the system through an external source while a mix of labeled and unlabeled \(X_2\) can return to \(X_1\), leave the system, or be converted into the following metabolite. This process continues for each metabolite in the pathway until \(X_N\). Unlabeled \(X_N\) enters the system through an external source while a mix of labeled and unlabeled \(X_N\) can return to the previous metabolite or leave the scope of the diagram}
	   \label{fig:Reversible Chain Reaction}
    \end{figure}

Consider the reversible chain reaction network in \cref{fig:Reversible Chain Reaction}. The scaled ODE model 
can be written in form \cref{re-scaled ODE} with the following variables and parameters. For \(n = 1,\ldots,N\),
		\begin{align*}
			\bar{x}_n &= \frac{X_n^\mathcal{U}}{X_n^\mathcal{T}}  &k_n &= \frac{F_{n}}{X_n^\mathcal{T}} &\alpha_n &= \frac{f_{3n-1}}{F_{n}}
		\end{align*}
and for \(n=1,\ldots, N-1\),
	\[\beta_{n,n+1} = \frac{f_{-(3n+1)}}{F_{n}}.\]

\subsubsection{Nondimensionalized Reversible Chain Reaction Model}\label{sec4.2.1}

To investigate the effectiveness of experimental data from a chain reaction network containing a disparity between time scales, we include the full nondimensionalization of the reversible chain reaction model to fast and slow subsystems. 
To simplify the analysis we will only consider the case where the turnover rate of one metabolite is much faster than the turnover rate of all other metabolites. 

Since each metabolite receives inputs only from neighboring metabolites, The matrix $B-I$ is tri-diagonal.
We will therefore construct the fast and slow systems using the system of equations in \cref{chain reaction system if eqns} instead of the matrix representation. 

    \begin{equation}\label{chain reaction system if eqns}
        \begin{aligned}
            \dot{\bar{x}}_1 &= k_1(-\bar{x}_1 +\beta_{1,2}\bar{x}_2+\alpha_1) \\
            \dot{\bar{x}}_i &= k_i((1-\alpha_i-\beta_{i,i+1})\bar{x}_{i-1}-\bar{x}_i+\beta_{i,i+1}\bar{x}_{i+1}+\alpha_i)\quad\text{for }i=2,\ldots,N-1\\
		      \dot{\bar{x}}_N  &= k_N((1-\alpha_N)\bar{x}_{N-1}-\bar{x}_N+\alpha_N)
        \end{aligned}
    \end{equation}
    
Assume we know the turnover rate of one metabolite in the reversible chain reaction network is much greater than the turnover rate of all other metabolites. 
Thus, we will have one fast variable and \(N-1\) slow variables. 
The effectiveness of the experimental data is directly related to the location of the metabolite with a quick turnover rate. 
We will derive the dimensionless model for two cases. 
In the first case, the metabolite with the fastest turnover rate is the first metabolite in the chain reaction network. 
In the second case, the metabolite with the fastest turnover rate is not the first in the chain reaction network.

\subsubsection*{Case 1:}\label{sec4.2.2}

Suppose the metabolite with the greatest turnover rate is the metabolite generated directly from nutrient, \(X_1\), making \(\bar{x}_1\) the fast variable. 
So, \(k_1 \gg k_i\) for \(i = 2,\ldots, N\). 
Let the dimensionless variable \(\tau = k_1t\) be the fast time scale and let \(\epsilon = k_N/k_1 \ll 1\). 
By taking \(\lim_{\epsilon\to 0}\), we obtain the fast subsystem.

    \begin{equation}\label{Fast Subsystem 1}
        \begin{aligned}
            \frac{d\bar{x}_1}{d\tau} &= -\bar{x}_1 +\beta_{1,2}+\alpha_1 \\
            \bar{x}_i &= 1 \quad\text{for }i=2,\ldots,N
        \end{aligned}
    \end{equation}
On the fast time scale, \(\bar{x}_1\) approaches the equilibrium of \cref{Fast Subsystem 1} at \(\bar{x}_1 = \alpha_1+\beta_{1,2}\). 

Let the dimensionless variable \(\eta = \epsilon \tau\) be the slow time scale. 
The slow system is derived by taking \(\lim_{\epsilon\to 0}\) to obtain 

    \begin{equation}\label{Slow Subsystem 1}
        \begin{aligned}
            0 &= -\bar{x}_1 +\beta_{1,2}\bar{x}_2+\alpha_1 \\
            \frac{d\bar{x}_i}{d\eta} &= \frac{k_i}{k_N}((1-\alpha_i-\beta_{i,i+1})\bar{x}_{i-1}-\bar{x}_i+\beta_{i,i+1}\bar{x}_{i+1}+\alpha_i)\quad\text{for }i=2,\ldots,N-1\\
		      \frac{d\bar{x}_N}{d\eta} &= (1-\alpha_N)\bar{x}_{N-1}-\bar{x}_N+\alpha_N
        \end{aligned}
    \end{equation}

In this case, the solution will start at the initial condition \(\vec{\boldsymbol{1}}\), then quickly approach the slow manifold, \(-\bar{x}_1 +\beta_{1,2}\bar{x}_2+\alpha_1=0\),  which is the nullcline of \(\bar{x}_1\). 
The quick jump from the initial condition to the slow manifold shows the fast dynamics of the system. 
Data collected quickly after the isotope switch will be on the fast time scale and may be used to estimate the fast turnover rate. 
Data collected later in the experiment is on the slow time scale \cref{Slow Subsystem 1} and may be used to approximate the other turnover rate parameters. 
The proportional flux parameters $\alpha_1$ and $\beta_{1,2}$ will be particularly easy to estimate because the solutions will remain on the slow manifold for most of the experiment. 
For this scenario, the fluxes approximated from data collected both early and later in the experiment will be accurate since both fast and slow dynamics can be observed.

\subsubsection*{Case 2:}\label{sec4.2.3}

On the other hand, if the metabolite with the greatest turnover rate is a downstream product of the metabolite generated directly from the nutrient, then fast dynamics will not be observed in the data. 
The issue arises because the initial condition is located on the slow manifold. 
Let \(X_{n}\) for \(n\neq 1\) be the metabolite with the greatest turnover rate, making \(\bar{x}_n\) the fast variable.
So, \(k_n\gg k_i\) for \(i = 1,\ldots, N\) with \(i\neq n\).
Let the dimensionless variable \(\tau = k_nt\) be the fast time scale and let \(\epsilon = k_1/k_n \ll 1\). 
By taking \(\lim_{\epsilon \to 0}\) we obtain the fast subsystem

    \begin{equation}\label{Fast Subsystem 2}
        \begin{aligned}
            \frac{d\bar{x}_n}{d\tau} &= 1-\bar{x}_n\\
            \bar{x}_i &= 1\quad\text{for }i=1,\ldots,N, i \neq n\\
        \end{aligned}
    \end{equation}

Given that the initial condition is \(\vec{1}\), the fast subsystem begins at its equilibrium \(\bar{x}_n=1\). 
This implies that the trajectory \((\bar{x}_1,\ldots,\bar{x}_N)\) does not leave the initial condition on the fast scale. 
Let the dimensionless variable \(\eta = \epsilon \tau\)  be the slow time scale. We derive the slow subsystem for this scenario by taking \(\lim_{\epsilon\to 0}\) 

    \begin{equation}\label{Slow Subsystem 2}
        \begin{aligned}
            \frac{d\bar{x}_1}{d\eta} &= -\bar{x}_1 +\beta_{1,2}\bar{x}_2+\alpha_1 \\
            0 &= (1-\alpha_n-\beta_{n,n+1})\bar{x}_{n-1}-\bar{x}_n+\beta_{n,n+1}\bar{x}_{n+1}+\alpha_n\\
            \frac{d\bar{x}_i}{d\eta} &= \frac{k_i}{k_1}(1-\alpha_i-\beta_{i,i+1})\bar{x}_{i-1}-\bar{x}_i+\beta_{i,i+1}\bar{x}_{i+1}+\alpha_i\quad\text{for }i=2,\ldots,N-1, i \neq n\\
		      \frac{d\bar{x}_N}{d\eta} &= \frac{k_N}{k_1}(1-\alpha_N)\bar{x}_{N-1}-\bar{x}_N+\alpha_N
        \end{aligned}
    \end{equation}

From \cref{Slow Subsystem 2} we obtain the slow manifold defined by \((1-\alpha_n-\beta_{n,n+1})\bar{x}_{n-1}-\bar{x}_n+\beta_{n,n+1}\bar{x}_{n+1}+\alpha_n = 0\) which is the \(\bar{x}_n\) nullcline. 
We know the slow manifold contains the initial condition by investigating the fast subsystem in \cref{Fast Subsystem 2}. 
Therefore solutions will display only slow dynamics and all variables will approach the equilibrium point at the intersection of all nullclines along the slow manifold. 
The turnover rate for the fast metabolite cannot be estimated. However, the proportional flux parameters for fluxes into the fast metabolite, $\alpha_n$ and $\beta_{n,n+1}$ will be evident throughout the experiment as the solutions track along the slow manifold. 
While the slow dynamics of the upstream metabolites limit our ability to estimate the fast dynamics accurately, the slow dynamics are still able to be estimated. 

The conclusions of case 1 and case 2 also apply to the irreversible chain reaction. In the irreversible chain reaction \(\beta_{i,i+1} = 0\). The initial condition is still on the slow manifold if a slow turnover rate is upstream of the faster turnover rate. The irreversible two-metabolite example illustrating this conclusion is included in the supplementary material.
    
\subsection{Two-Metabolite Example}\label{sec4.2.4}

Using the reversible two-metabolite example, we will illustrate both cases in \cref{sec4.2.1}. Consider the reversible two-metabolite reaction network in \cref{fig:Reversible 2 metabolites} and its scaled ODE model \cref{re-scaled R2M ODE}. 

\subsubsection*{Case 1: Fast-Slow System}
Suppose we know the turnover rate of \(X_1\) is much greater than the turnover rate of \(X_2\), i.e. \(k_1 \gg k_2\). 
Let the dimensionless variable \(\tau = k_1t\) be the fast time scale and let \(\epsilon = k_2/k_1 \ll 1\). Following the same process as in \cref{sec4.2.1} for the general reversible chain reaction model, we obtain the fast subsystem.
	\begin{align*}
  \frac{d\bar{x}_1}{d\tau} &= -\bar{x}_1 + \beta_{2,1} + \alpha_1\\
  \bar{x}_2 &= 1
	\end{align*}
Now, let \(\eta = \epsilon\tau\) be the slow time scale to obtain the slow subsystem.
	\begin{align*} 
  0 &= -\bar{x}_1 + \beta_{1,2} \bar{x}_2 + \alpha_1 \\
  \frac{d\bar{x}_2}{d\eta} &= (1-\alpha_2)\bar{x}_1 -\bar{x}_2 + \alpha_2
	\end{align*}
In this example, the slow manifold is defined by \(\bar{x}_2 = \frac{1}{\beta_{1,2}}(\bar{x}_1-\alpha_1)\).

    \begin{figure}[H]
	   \centering
	   \includegraphics[width=0.65\textwidth]{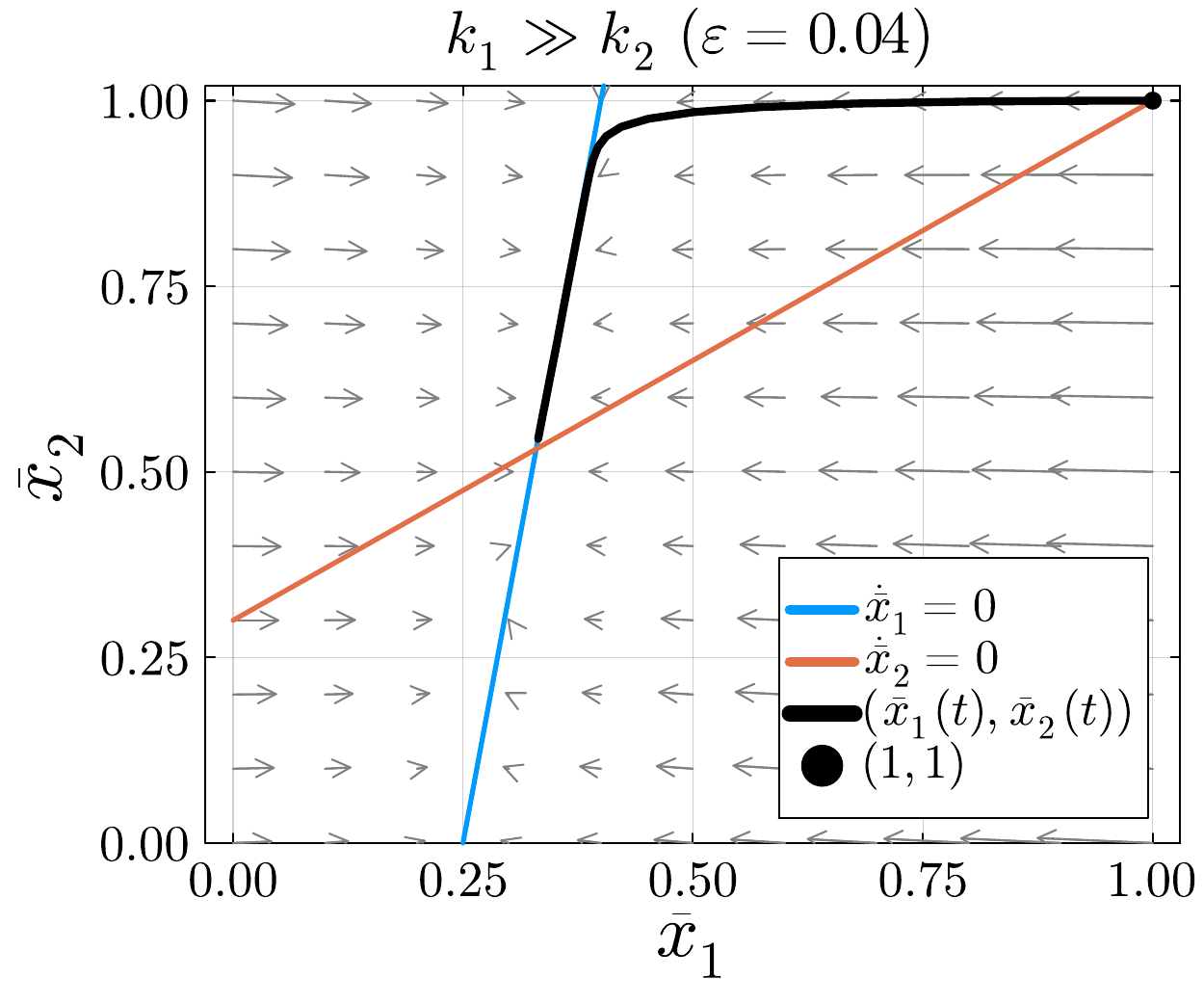}
	   \caption{The phase plane for the reversible two-metabolite model in the case where the turnover rate of \(X_1\) is much greater than the turnover rate of \(X_2\) when including the reversible reaction. The parameter values used to create this figure are \(k_1= 1\), \(\alpha_1 = 1/4\), \(\beta_{1,2}=3/20\), \(k_2= 1/25\), and \(\alpha_2= 3/10\)}
	   \label{fig:Reversible FS phase plane}
    \end{figure}

The solution \((\bar{x}_1(t),\bar{x}_2(t))\) in \cref{fig:Reversible FS phase plane} begins at the initial condition \((1,1)\). Next, it quickly moves from the initial condition to the slow manifold. This portion of the solution shows the fast dynamics of the system. After approaching the slow manifold, the solution slowly approaches the equilibrium point \cref{irreversible SS}. 
The movement along the slow manifold represents the slow dynamics of the system. 
In this case, data collected quickly after the isotope switch will be on the fast time scale and data collected later in the experiment will be on the slow time scale. 
With both early and later measurements from this reaction network, KFP will provide more accurate flux estimations given that both fast and slow dynamics are represented.

\begin{figure}[H]
        \centering
        \includegraphics[width=0.9\textwidth]{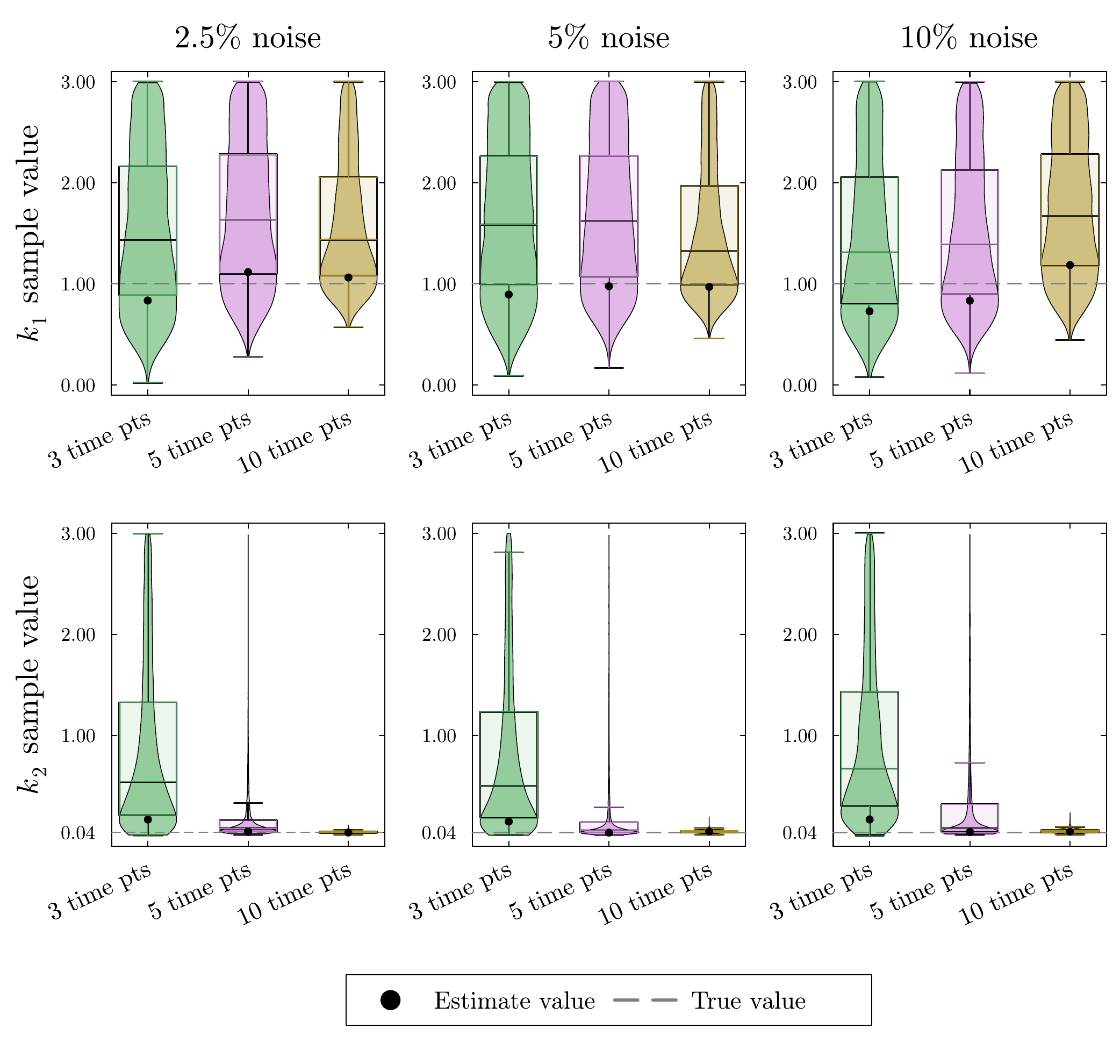}
        \caption{Violin and box plots for posterior distributions of \(k_1\) and \(k_2\) for the reversible two-metabolite fast-slow system. The simulated data are taken from solutions of the system with parameters \(k_1 = 1\), \(\alpha_1 = 1/4\), \(\beta_{1,2} = 3/20\), \(k_2 = 1/25\), \(\alpha_2 =3/10\) at 3, 5 or 10 equally distributed time points. At each time point, three noisy data points are simulated by adding random normally distributed noise to the true solutions with a standard deviation of 2.5, 5, or 10\% of the true value. Since \(k_1\), and \(k_2\) are turnover rates, we use the naive prior distributions
		$k_1 \sim U(0,3),\; k_2 \sim U(0,3)$. The mode of the posterior samples is taken as the estimated value. Outliers are excluded in the box plots}
        \label{fig:R2MFastSlowDenK1K2}
    \end{figure}

We simulate data for this system as described in Appendix \cref{appendix: Data and BPE}. The resulting posterior sample distributions for \(k_1\) and \(k_2\) are included in \cref{fig:R2MFastSlowDenK1K2}. As with previous Bayesian estimates, all parameters including proportional flux parameters and turnover rates are given naive priors. In this section, we show only the posterior distributions for the turnover rates since the estimation of the proportional flux parameters follows the pattern seen in \cref{sec3}. The simulated data, the estimates for the proportional flux parameters, and the credible intervals for each parameter estimation can all be found in the supplement. 

In the top row of \cref{fig:R2MFastSlowDenK1K2}, we see the estimates for \(k_1\) are accurate but highly uncertain. 
Because there was only a single time point taken during the rapid decay of \(\bar{x}_1\), we cannot get a good estimate for \(k_1\). 
In the bottom row of \cref{fig:R2MFastSlowDenK1K2}, we see the estimate for \(k_2\) increases in accuracy and decreases in uncertainty with the addition of time points despite the noise level. 

\subsubsection*{Case 2: Slow-Fast System}
Now, we assume the opposite situation. 
Suppose we know the turnover rate of \(X_2\) is much greater than the turnover rate of \(X_1\), i.e. \(k_2 \gg k_1\). 
Let the dimensionless variable \(\tau = k_2t\) be the fast time scale and let \(\epsilon = k_1/k_2 \ll 1\). 
Following the same process as in \cref{sec4.2.1} we obtain the fast subsystem.

	\begin{align*} 
  \bar{x}_1 &= 1\\
  \frac{d\bar{x}_2}{d\tau} &= 1 - \bar{x}_2
	\end{align*}
 
Now, let \(\eta = \epsilon\tau\) be the slow time scale to obtain the slow subsystem.

	\begin{align*}
  \frac{d\bar{x}_1}{d\eta} &= -\bar{x}_1 + \beta_{1,2} \bar{x}_2 + \alpha_1 \\
  0 &= (1-\alpha_2)\bar{x}_1 - \bar{x}_2 + \alpha_2
	\end{align*}
In this example, the slow manifold is defined by the line \(\bar{x}_2 = (1-\alpha_2)\bar{x}_1+\alpha_2\).

    \begin{figure}[H]
	   \centering
	   \includegraphics[width=0.65\textwidth]{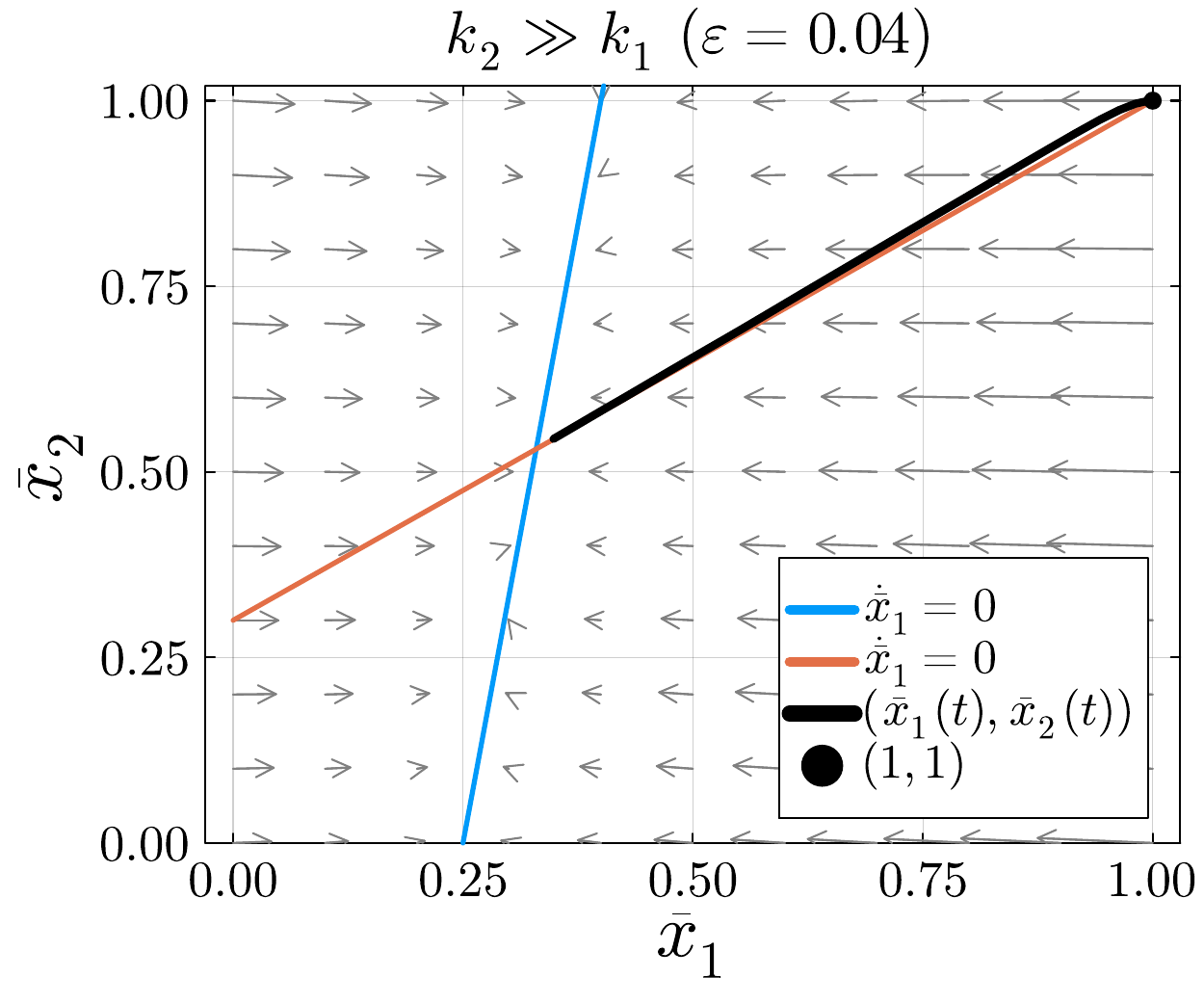}
	   \caption{The phase plane for the reversible two-metabolite model in the case where the turnover rate of \(X_2\) is much greater than the turnover rate of \(X_1\) when including the reversible reaction. The parameter values used to create this figure are \(k_1= 1/25\), \(\alpha_1 = 1/4\), \(\beta_{1,2}=3/20\), \(k_2= 1\), and \(\alpha_2= 3/10\)}
	   \label{fig:Reversible SF phase plane}
    \end{figure}

The solution \((\bar{x}_1(t), \bar{x}_2(t))\) in \cref{fig:Reversible SF phase plane} begins at the initial condition \((1,1)\). 
The initial condition is located on the slow manifold \(\bar{x}_2 = (1-\alpha_2)\bar{x}_1 + \alpha_2\). 
The solution will slowly move along the slow manifold until it approaches the equilibrium located at the intersection of \(\dot{\bar{x}}_1 = 0\) and \(\dot{\bar{x}}_2 = 0\). 
Therefore, only the slow dynamics will be visible in the experimental data. 
Metabolite \(X_1\)'s slow turnover rate will restrict the amount of label entering metabolite \(X_2\), preventing an accurate approximation of the much faster rate through \(X_2\). 
To measure the faster turnover rate through \(X_2\), we would need to change the initial condition of the system. 
Changing the initial condition would require a different experimental setup.

\begin{figure}[H]
        \centering
        \includegraphics[width=0.9\textwidth]{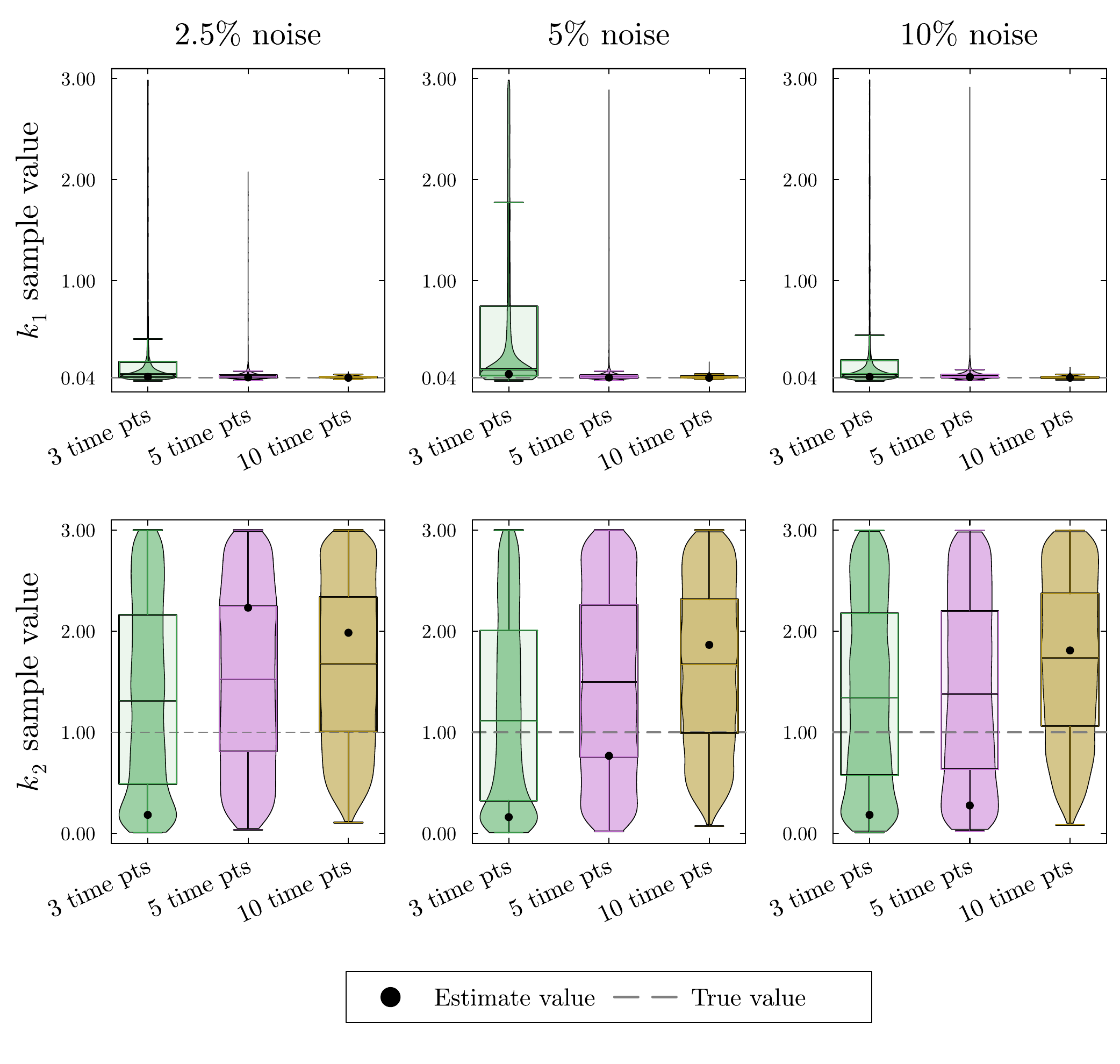}
        \caption{Violin and box plots for posterior distributions of \(k_1\) and \(k_2\) for the reversible two-metabolite slow-fast system. The simulated data are taken from solutions of the system with parameters \(k_2 = 1/25\), \(\alpha_1 = 1/4\), \(\beta_{1,2} = 3/20\), \(k_2 = 1\), \(\alpha_2 =3/10\) at 3, 5 or 10 equally distributed time points. At each time point, three noisy data points are simulated by adding random normally distributed noise to the true solutions with a standard deviation of 2.5, 5, or 10\% of the true value. Since \(k_1\), and \(k_2\) are turnover rates, we use the naive prior distributions
		$k_1 \sim U(0,3),\; k_2 \sim U(0,3)$. The mode of the posterior samples is taken as the estimated value. Outliers are excluded in the box plots}
        \label{fig:R2MSlowFastDenK1K2}
    \end{figure}

We simulate data for this system as described in Appendix A. The resulting posterior sample distributions for \(k_1\) and \(k_2\) are included in \cref{fig:R2MSlowFastDenK1K2}.

At the top of \cref{fig:R2MSlowFastDenK1K2}, the estimates for \(k_1\) are very accurate and have extremely low uncertainty. 
With many experimental measurements throughout the decay of \(\bar{x}_1\), we expected to get good estimates for \(k_1\). 
Contrarily, we receive very poor estimates for \(k_2\) with incredibly high uncertainty as shown at the bottom of \cref{fig:R2MSlowFastDenK1K2}. 
The solution for \(\bar{x}_2\) is not sensitive to changes in the large \(k_2\) value. 
With enough data points, Bayesian parameter estimation can rule out small values in \(k_2\) but cannot define a good range for the true value of the parameter.

\section{Discussion}\label{sec5}

In this paper, we have provided a general formulation of KFP models along with a useful scaling that eliminates dependence on total metabolite concentrations. 
This general formulation allows for more complex graph structures and facilitates the understanding of the relationship between graph structure and parameter identifiability. We have shown that steady-state data alone can only provide estimates of relative fluxes, not turnover rates, and only when the steady-state system is not underdetermined. 
Finally, by investigating a separation of time scales, we show that the fast turnover rate can only be estimated if it is the node that receives the labeled input. 
In these latter sections, we have demonstrated the utility of Bayesian parameter estimation in conjunction with KFP. 
In contrast to other fitting methods, Bayesian methods allow the experimenter to judge not only the best fitting parameter but also to quantify the credibility of that estimate. 
In addition, by identifying larger turnover rates with large credible intervals, fast turnover metabolites can be identified and targeted for future experiments.
We answered three crucial questions for the experimenter, but many open questions have yet to be investigated. 
For the remainder of this section, we will address a few of the many open questions related to the theme of how to incorporate multiple experiments or multiple types of data. In some cases, it may be reasonable to use several types of data in one parameter estimation. In others, the Bayesian approach allows for an iterative method where new experiments are used to update the distribution for each parameter value. 
 
Often it is not feasible to accurately measure all metabolites in the pathway of interest with the mass spectrometer. 
Many unknown and uncontrollable features factor into the ability to detect a metabolite. For example, some metabolites are very low in abundance, others ionize poorly and still others are unstable, all of which can prevent identification in the sample. 
Therefore, it is important to consider how missing measurements affect the identifiability of the model's parameters. 
If the unmeasurable metabolite is included in the network and therefore in the KFP model, the steady-state problem becomes underdetermined. 
By excluding the unmeasurable metabolite from the KFP model we reduce the number of parameters. 
With a reduced model, the experimenter must be aware that the reduced parameter set is a combination of some of the unreduced parameters. 
Therefore, an important open question is how the estimates in the reduced model relate to the parameters in the unreduced model.

All of the examples presented here have only one source of labeled input (i.e. \(|\mathcal{E}_{\mathcal{L}}| = 1)\)). 
In practice, the generation of secondary tracers can lead to more than one source for the label to enter the system simultaneously. 
For a more biologically accurate representation, multiple sources of labeled input should be included in the model. In the general model formulation we can incorporate multiple labeled inputs. However, 
we suspect this will further complicate the parameter estimation problem. 
For example, if the type of label from one source does not differ from the other we will not be certain which source provided the label. 

Currently, the KFP models only recognize one type of label. 
The metabolites are either labeled or unlabeled. 
In an isotope tracing experiment, the total pool of a metabolite can be composed of multiple types of label. 
If a metabolite contains multiple carbon atoms, those carbons may be a mix of \(^{13}\)C and unlabeled \(^{12}\)C atoms. For example, citrate contains 6 carbon atoms so there are 7 different types of label called isotopologues for that metabolite (i.e. 0, 1, 2, 3, 4, 5, or 6 \(^{13}\)C molecules). 
If the labeled substrate is pyruvate, one sequence of reactions including pyruvate dehydrogenase can produce citrate with 2 \(^{13}\)C atoms while a different route through pyruvate carboxylase produces citrate with 3 \(^{13}\)C atoms. 
To more accurately represent the biological system, multiple types of labels, should be incorporated into the KFP model. Specific labeling types can be incorporated by including a variable for each label type for each metabolite. This greatly increases the dimension of the KFP model, but using the stoichiometry of each reaction, may lead to better flux estimates.

Another example where multiple labels may be useful is to include experiments with different stable isotopes. The same metabolite could have labeled nitrogen, oxygen or hydrogen atoms. By repeating the tracing experiment with the same pathway but with multiple isotopes we can better deduce the sequence of reactions for each molecule. 
Alternatively, one may want to use data from multiple experiments studying the same pathway but differing in the source node of the labeled input. 
This would create an ensemble of KFP models working together to estimate the same set of unknown fluxes. 
The ensemble may produce more accurate estimates because the fluxes estimated by one model can be used as to inform the choice of prior distributions for the following model. 
Iterating through the models may help the parameter estimation method reduce uncertainty and pinpoint the true parameter values. 

Previously, KFP has been used to estimate relative flux changes (rKFP) \citep{huang2014estimating} between an experimental and control condition. rKPF requires multiple parameter fits on the same KFP model. 
It would be useful to consider how the conditions and constraints described here will impact the comparison of multiple experimental conditions with the same pathway. 

Kinetic Flux Profiling is an intriguing application for parameter identifiability because the accurate measurement of the flux parameters is the central motivation for the method. 
The algebraic constraints on the parameters provide additional structure to the problem that makes the parameter estimation not only more tractable but also more nuanced. 
From a biological standpoint, the information obtained from the estimates of flux can provide essential insight into the mechanistic changes in metabolism during disease. 
At the same time, the limited and noisy data challenge our ability to make accurate estimates. 

\backmatter

\bmhead{Supplementary information}
Supplementary file 1 includes additional figures and examples resulting from Bayesian parameter estimation.
\href{run:./supplementary_file_1.pdf}{Supplementary File 1}

\bmhead{Acknowledgements}
We would like to thank Mitchell Riley for helping review and improve the manuscript.

\section*{Declarations}

The authors have no relevant financial or non-financial interests to disclose.

The authors declare that the data supporting the findings of this study are available within the paper or its supplementary information files.

This work was supported by grants NIH RO1 DK104998 (EBT), NIH RO1 DK138664 (EBT), University of Iowa Distinguished Scholar Award (EBT)

\begin{appendices}

\section{Simulated Data and Bayesian Parameter Estimation}\label{appendix: Data and BPE}

Using Bayesian parameter estimation, we illustrate the accuracy and reliability of parameter estimation for various scenarios. 
To achieve this, we vary the number of time points and level of noise in the simulated experimental data for each model.

First, we select parameter values to compute the true solution and simulate experimental measurements by adding noise to the true solution at specified time points. 
We vary the number of measurements using 3, 5, and 10 time points evenly spaced throughout the time frame of the experiment. 
Next, we add 2.5\%, 5\%, or 10\% random normally distributed noise to the true solution at these time points. 
Since an experiment is usually run three times, we have 3 noisy measurements for each time point for each metabolite measured. 
We run Bayesian estimation using measurements from all three experiments instead of the average value of the measurements from experiments first.

We apply Bayesian parameter estimation to our models using the Julia package \verb|Turing.jl| \citep{ge2018t}. 
In the process of parameter estimation, we will need to solve the differential equation many times. 
We will be using the Julia package \verb|DifferentialEquations.jl| \citep{DifferentialEquations.jl-2017} for defining and solving the DE.

Bayesian parameter estimation requires prior distributions for each parameter. 
Since we expect the experimenter to have little prior knowledge of the true parameter values, we provide the method with naive prior distributions. 
For \(k_1\) and \(k_2\), we set the prior distributions to be
		\[k_X \sim U(0,3) \quad\text{and}\quad k_Y\sim U(0,3). \] 
	Since \(\alpha_1\), \(\alpha_2\), and \(\beta_{2,1}\) are proportions we naturally set their prior distributions to be
		\[\alpha_1 \sim U(0,1),\quad\alpha_2 \sim U(0,1), \quad\text{and}\quad \beta_{2,1} \sim U(0,1).\]
  
This method returns samples from the posterior distributions for each parameter. 
We use the No-U-Turn sampler (NUTS) to sample the posterior distribution. 
NUTS is a Markov Chain Monte Carlo (MCMC) algorithm that builds a set of points that covers a wide range of the desired distribution by automatically adaption the step size and number of steps \citep{hoffman2014no}. 
After receiving sample distributions for our parameter values, we must choose a metric from the distribution to be the parameter estimate. 
We use the mode of distribution to the estimated value. The mode is computed by finding the maximum of the kernel density estimation of the sample distribution. 
Uncertainty in the estimated value is measured by the spread of the distribution. 
If the difference in the 97.5\% and 2.5\% quartiles is large, we will say the estimate has high uncertainty. 
Conversely, if the difference between the 97.5\% and 2.5\% quartiles is small the estimate has low uncertainty. The difference between the 97.5\% and 2.5\% quartile is called the Bayesian 95\% credible interval. 
The 95\% credible interval is the interval with a 95\% probability of containing the true parameter value given the data \citep{hespanhol2019understanding}.

\end{appendices}

\bibliography{sn-bibliography}
\end{document}